\documentclass[12pt,a4paper]{article}
\usepackage[T2A]{fontenc}
\usepackage[cp1251]{inputenc}
\usepackage[english,russian]{babel}
\usepackage[tbtags]{amsmath}
\usepackage{amsfonts,amssymb,mathrsfs,amscd}

\usepackage{srcltx}\usepackage{graphicx,color%,ifsym%,txfonts
}
\usepackage[all]{xy}
\usepackage{latexsym,amsfonts,amssymb,amsmath,longtable,%amsthm,
mathrsfs,theorem,cite
}
\usepackage{stmaryrd}
\usepackage[hyper]{amsbib}

\numberwithin{equation}{section}
\theoremstyle{plain}
\newtheorem{theorem}{Теорема}
\newtheorem{lemma}{Лемма}[section]
\newtheorem{propos}{Предложение}
\newtheorem{corollary}{Следствие}
\def\h{\operatorname{h}}
\def\m{\operatorname{m}}
\def\sm{\operatorname{sm}}
\def\Hall{\operatorname{Hall}}
\def\SL{\operatorname{SL}}
\def\Sp{\operatorname{Sp}}
\def\GL{\operatorname{GL}}
\def\PSL{\operatorname{PSL}}
\def\PGL{\operatorname{PGL}}
\def\PSp{\operatorname{PSp}}
\def\Sym{\operatorname{Sym}}
\def\Alt{\operatorname{Alt}}
\def\Aut{\operatorname{Aut}}
\def\Inn{\operatorname{Inn}}
\def\Syl{\operatorname{Syl}}
\def\qed{\hfill $\square$}
\def\E{\mathscr{E}}
\def\C{\mathscr{C}}
\def\D{\mathscr{D}}
\def\M{\mathscr{M}}
\begin{document}

%\title{\vspace{-1.5cm}\\ 
\begin{center}
{\Large\bf When is the search of relatively maximal subgroups reduced to quotients?
\footnote{W.~Guo was supported by the National Natural Science Foundation of China (No.12171126), Wu Wen-Tsun Key Laboratory of Mathematics of Chinese Academy of Sciences and Key Laboratory of Engineering Modeling and Statistical Computation of Hainan Province. D.\,O.~Revin was supported by RFBR and BRFBR, project  No.~20-51-00007 and by the Program of Fundamental Scientific Research of the SB RAS~I.1.1., project~No.~0314-2016-0001.}\\
}%Когда поиск относительно максимальных подгрупп редуцируется к факторгруппам?
%}\\
%\author

~

~

{\bf Wenbin Guo
and
Danila O.~Revin}
%\date{}
\end{center}

\def\abstractname{\bf Abstract}

%\maketitle

\begin{abstract}
Let ${\mathfrak{X}}$ be a class of finite groups closed under taking subgroups, homomorphic images, and extensions. Denote by ${\mathrm{k}}_{\mathfrak{X}}(G)$ the number of conjugacy classes ${\mathfrak{X}}$-maximal subgroups of a finite group $G$. The natural problem to describe up to conjugacy ${\mathfrak{X}}$-maximal subgroups of a given finite group is complicated by the fact that it is not inductive. In particular, generally speaking, the image of an ${\mathfrak{X}}$-maximal subgroup is not ${\mathfrak{X}}$-maximal in the image of a homomorphism. Nevertheless, there are group homomorphisms which preserve the number of conjugacy classes of ${\mathfrak{X}}$-maximal subgroups (for example, the homomorphisms whose kernels are ${\mathfrak{X}}$-groups). Under such homomorphisms, the image of an ${\mathfrak{X}}$-maximal subgroup is always ${\mathfrak{X}}$-maximal and, moreover, there is a natural bijection between the conjugacy classes of ${\mathfrak{X}}$-maximal subgroups of the image and preimage. All such homomorphisms are completely described in the paper. More precisely, it is proved that, for a homomorphism $\phi$ from a group $G$, the equality ${\mathrm{k}}_{\mathfrak{X}}(G)={\mathrm{k}}_{\mathfrak{X}}(\mathrm{im}\, \phi)$ holds if and only if ${\mathrm{k}}_{\mathfrak{X}}(\ker \phi)=1$, which in turn is equivalent to the fact that the composition factors of the kernel of $\phi$ belong to an explicitly given list.

\smallskip

\end{abstract}

\section{Введение}\label{Intro}

\subsection{Основной результат}\label{MainRes}

В работе рассматриваются только конечные группы, и термин ``группа'' используется в значении ``конечная группа''.

Группа, принадлежащая классу групп~${\mathfrak{X}}$, называется {\it ${\mathfrak{X}}$-группой}. Множество {\it максимальных} (по включению) {\it ${\mathfrak{X}}$-подгрупп} (иначе {\it ${\mathfrak{X}}$-мак\-си\-маль\-ных подгрупп}) группы~$G$ обозначается через $\m_{\mathfrak{X}}(G)$. Сама группа~$G$, действуя на $\m_{\mathfrak{X}}(G)$ сопряжениями, разбивает это множество на орбиты~--- классы сопряженности. Число этих классов обозначим через~${\mathrm{k}}_{\mathfrak{X}}(G)$. Термин ``{\it относительно максимальная подгруппа}'', использованный в названии статьи, предложен Х.\,Виландом~\cite{WieCanb} для того, чтобы отличать ${\mathfrak{X}}$-максимальные для некоторого класса ${\mathfrak{X}}$ подгруппы от просто максимальных (т.\,е. подгрупп, максимальных среди собственных).  Следуя Виланду \cite{Wie4,Wie3}, мы говорим, что непустой класс ${\mathfrak{X}}$ конечных групп является {\it полным}, если он замкнут относительно взятия подгрупп, гомоморфных образов и расширений. Последнее означает, что ${G\in{\mathfrak{X}}}$, как только группа $G$ содержит такую нормальную подгруппу $N$, что $N$ и $G/N$~--- ${\mathfrak{X}}$-группы.

В случае полного класса ${\mathfrak{X}}$ вопрос о том, {\it когда поиск ${\mathfrak{X}}$-максимальных подгрупп группы $G$ можно редуцировать к аналогичной задаче для факторгруппы $G/N$}, сформулированный в названии работы, оказывается эквивалентным вопросу {\it когда ${\mathrm{k}}_{\mathfrak{X}}(G)={\mathrm{k}}_{\mathfrak{X}}(G/N).$} Ответить на него позволяет

\begin{theorem}\label{Reduktionssatz_main}
Пусть ${\mathfrak{X}}$~--- полный класс групп и $N$~--- нормальная подгруппа конечной группы $G$. Тогда $$\quad{\mathrm{k}}_{\mathfrak{X}}(G)={\mathrm{k}}_{\mathfrak{X}}(G/N)\quad\text{если и только если}\quad{\mathrm{k}}_{\mathfrak{X}}(N)=1.$$
\end{theorem}

Другими словами, число классов сопряженных ${\mathfrak{X}}$-максимальных подгрупп тогда и только тогда не меняется при переходе от группы $G$ к факторгруппе $G/N$, когда в~$N$ все ${\mathfrak{X}}$-максимальные подгруппы сопряжены.

Утверждение ``если'' доказано в~\cite{GRV_Reducktionssatz}. Там же получено исчерпывающее описание групп $A$ с ${\mathrm{k}}_{\mathfrak{X}}(A)=1$: это условие равносильно тому, что каждый неабелев композиционный фактор группы $A$ либо принадлежит ${\mathfrak{X}}$, либо изоморфен простой группе, явным образом указанной в~\cite[приложение~A]{GRV_Reducktionssatz}. Тем самым, в случае полного класса ${\mathfrak{X}}$ на вопрос в названии работы удается дать исчерпывающий ответ.

Предположение о полноте класса ${\mathfrak{X}}$ в теореме~\ref{Reduktionssatz_main} существенно. Теорема неверна, если ${\mathfrak{X}}=\mathfrak{A}$~--- класс всех абелевых или ${\mathfrak{X}}=\mathfrak{N}$~--- класс всех нильпотентных групп. В обоих случаях ${\mathfrak{X}}$ не замкнут относительно расширений и   $${{\mathrm{k}}_{\mathfrak{X}}(\Sym_3)=2}\,\,\ne\,\,{1={\mathrm{k}}_{\mathfrak{X}}(\Sym_3/\Alt_3)},\quad \text{хотя} \quad {\mathrm{k}}_{\mathfrak{X}}(\Alt_3)=1.$$

\subsection{Мотивировка и предыстория}\label{Motivation}

Начиная с основополагающих работ Э.\,Галуа \cite{Galois} и К.\,Жордана  \cite{Jordan0,Jordan}, в теории конечных групп и ее приложениях мы сталкиваемся со следующим типом задач.
 \begin{itemize}
   \item {\it Дана группа $G$ (например, симметрическая) и класс ${\mathfrak{X}}$ конечных групп (например, класс разрешимых групп); требуется найти ${\mathfrak{X}}$-подгруппы группы $G$.}
 \end{itemize}
 В общей постановке решать задачи подобного типа вряд ли возможно. Когда, подобно классу разрешимых групп, класс ${\mathfrak{X}}$ является полным, можно ограничиться поиском  ${\mathfrak{X}}$-максимальных подгрупп.

В дальнейшем символ ${\mathfrak{X}}$ всегда используется  для обозначения некоторого фиксированного полного класса. Помимо тривиальных примеров, к которым относятся класс $\mathfrak{G}$ всех групп и класс $\mathfrak{E}$ групп порядка~1, а также помимо упомянутого выше класса $\mathfrak{S}$ разрешимых групп,  наиболее типичными примерами полных классов являются классы $\mathfrak{G}_\pi$ всех $\pi$-групп и  $\mathfrak{S}_\pi$ всех разрешимых $\pi$-групп для данного подмножества $\pi$ множества $\mathbb{P}$ всех простых чисел (напомним, что {\it $\pi$-группа}~--- это группа, у которой любой простой делитель порядка принадлежит~$\pi$). Эти два класса в определенном смысле экстремальны, поскольку для любого полного класса ${\mathfrak{X}}$ выполнены включения
$$
\mathfrak{S}_\pi\subseteq{\mathfrak{X}}\subseteq\mathfrak{G}_\pi,
$$
где через $\pi$ обозначено множество $$\pi({\mathfrak{X}})=\{p\in\mathbb{P}\mid \text{ существует } G\in{\mathfrak{X}}, \text{ для которой } p \text{ делит } |G|\}.$$
Полными являются также классы $\pi$-отделимых, $\pi$-разрешимых (в частности, $p$-раз\-ре\-ши\-мых) групп\footnote{Напомним, что  группа называется {$\pi$-отделимой}, если она обладает (суб)нормальным рядом, у которого все факторы являются $\pi$- или $\pi'$-группами, где $\pi'=\mathbb{P}\setminus\pi$. Если при этом все $\pi$-факторы этого ряда разрешимы, группа называется {\it $\pi$-разрешимой}. Особенно важную роль в теории групп играют $p$-разрешимые группы (случай, когда $\pi=\{p\}$), исследованные в классической работе Ф.\,Холла и Г.\,Хигмэна~\cite{HallHigman}. Эта работа позволила свести ослабленную проблему Бернсайда к группам примарного периода, а ее  методы легли в основу локального теоретико-группового анализа, см. работы Дж.\,Томпсона \cite{Th1,Th2}, его эссе, посвященное Холлу~\cite{ThompsonEssay}, и его знаменитую совместную работу с У.Фейтом о разрешимости групп нечетного порядка \cite{FeitThompsonOddOrder}.}.

 Естественно изучать ${\mathfrak{X}}$-максимальные подгруппы с точностью до сопряженности. Мы будем называть {\it ${\mathfrak{X}}$-схемой} группы $G$ полную систему представителей ее классов сопряженных ${\mathfrak{X}}$-максимальных подгрупп. Другими словами, ${\mathfrak{X}}$-схема группы --- это максимальная по включению система попарно несопряженных ${\mathfrak{X}}$-максимальных подгрупп. Мощность ${\mathfrak{X}}$-схемы группы $G$ равна  ${\mathrm{k}}_{\mathfrak{X}}(G)$~--- числу классов сопряженности ${\mathfrak{X}}$-максимальных подгрупп. Поиск ${\mathfrak{X}}$-схемы и описание строения ее элементов может рассматриваться как основная цель в задачах выделенного типа.

  Успешно решать такие задачи можно, если ${\mathfrak{X}}=\mathfrak{G}_p$ --- класс $p$-групп для некоторого простого числа~$p$. Тогда ${\mathfrak{X}}$-максимальные подгруппы оказываются силовскими $p$-подгруппами, и такие подгруппы в любой группе сопряжены~\cite{Sylow}. Также в разрешимых группах ${\mathfrak{X}}$-максимальные подгруппы --- это, в точности, $\pi({\mathfrak{X}})$-холловы подгруппы, которые образуют один класс сопряженности в соответствии с теоремой Ф.\,Холла \cite{Hall1}. Искать силовские и холловы подгруппы в произвольной конечной группе существенно помогают их хорошие свойства, позволяющие переходить от группы к секциями ее нормального или субнормального ряда и использовать соображения индукционного характера. Скажем, если $H$~--- силовская $p$-подгруппа группы $G$, то $H\cap N$ и $HN/N$~--- силовские $p$-подгруппы в $N$ и $G/N$ соответственно для любой $N\trianglelefteqslant G$.

К числу сложностей, возникающих в случае ${\mathfrak{X}}\ne\mathfrak{G}_p$, относится выраженно неиндуктивный характер рассматриваемых задач. Пересечение $H\cap N$ подгруппы $H\in\m_{\mathfrak{X}}(G)$ c нормальной подгруппой $N$ и образ $HN/N$ в факторгруппе $G/N$ (эквивалентно, образ $H$ при произвольном эпиморфизме) в общем случае уже не являются ${\mathfrak{X}}$-максимальными подгруппами в $N$ и $G/N$.

Виланд \cite{Wie3} заметил, что  для пересечений с нормальными и даже субнормальным подгруппами ситуацию можно отчасти исправить, если несколько расширить понятие ${\mathfrak{X}}$-мак\-си\-маль\-ной подгруппы и изучать т.\,н. ${\mathfrak{X}}$-суб\-мак\-си\-маль\-ные подгруппы\footnote{ Подгруппу $H$ группы $G$ Виланд назвал {\it ${\mathfrak{X}}$-суб\-мак\-си\-маль\-ной} (пишем $H\in\sm_{\mathfrak{X}}(G)$), если $G$ можно вложить в качестве субнормальной подгруппы в некоторую группу $G^*$ так, чтобы для подходящей $H^*\in\m_{\mathfrak{X}}(G^*)$ выполнялось равенство $H=H^*\cap G$. Из определения следует, что если $H\in\sm_{\mathfrak{X}}(G)$ и $N\trianglelefteqslant\trianglelefteqslant G$, то $H\cap N\in\sm_{\mathfrak{X}}(N)$. Понятие ${\mathfrak{X}}$-суб\-мак\-си\-маль\-ной подгруппы оказывается содержательным и отличным от понятия просто ${\mathfrak{X}}$-подгруппы благодаря теореме Виланда-Хартли \cite[5.4(a)]{Wie3}, \cite[теорема~2]{RSV}, утверждающей, что {\it если ${H\in\sm_{\mathfrak{X}}(G)}$, то~${H\in\m_{\mathfrak{X}}(N_G(H))}$}. }.

 В данной работе мы сфокусируем внимание на поведении ${\mathfrak{X}}$-максимальных подгрупп при гомоморфизмах.
 Попытка расширить понятие ${\mathfrak{X}}$-максимальной подгруппы так, чтобы оно было согласовано с гомоморфными образами, возвращает нас к необходимости изучать все ${\mathfrak{X}}$-подгруппы, см.~\cite{Epimax}.  Оказывается, что если для класса ${\mathfrak{X}}$ существует группа $L$ с несопряженными ${\mathfrak{X}}$-мак\-си\-маль\-ными подгруппами\footnote{Существование такой группы эквивалентно тому, что ${\mathfrak{X}}$ отличен от классов $\mathfrak{G}$, $\mathfrak{E}$ и $\mathfrak{G}_p$ для любого простого числа~$p$, см.~\cite{Epimax}.}, то любая группа $G$ является образом некоторого гомоморфизма (а~именно,  естественного эпиморфизма из регулярного сплетения~${L\wr G}$), при котором {\it всякая} (не только ${\mathfrak{X}}$-максимальная) ${\mathfrak{X}}$-подгруппа совпадает с образом некоторой ${\mathfrak{X}}$-максимальной подгруппы \cite[14.2]{Wie4}, \cite[4.3]{Wie3}.

Другая очевидная сложность, возникающая при переходе от группы к факторгруппе или эпиморфному образу состоит в том, что несопряженные ${\mathfrak{X}}$-мак\-си\-маль\-ные подгруппы    могут иметь сопряженные образы (например, при тривиальном эпиморфизме $G\rightarrow 1$). Тем самым может быть утрачена информация о классах сопряженности ${\mathfrak{X}}$-мак\-си\-маль\-ных подгрупп.

Поэтому важно выделить те случаи, когда переход от группы $G$ к факторгруппе $G/N$ является редукцией для выделенного типа задач: когда сохраняется ${\mathfrak{X}}$-мак\-си\-маль\-ность подгрупп и не искажается информация об их сопряженности. Другими словами, когда ${\mathfrak{X}}$-схема переходит в ${\mathfrak{X}}$-схему. Очевидным примером служит ситуация, когда $N\in{\mathfrak{X}}$. Менее очевидные случаи указал Виланд в~\cite{Wie4}. Там же он наметил программу, как
  найти все такие  ``хорошие'' случаи. Настоящей статьей мы завершаем исследования по этой программе.

Именно сохранение мощности ${\mathfrak{X}}$-схемы, т.\,е. равенство
\begin{equation}\label{1}
{\mathrm{k}}_{\mathfrak{X}}(G)={\mathrm{k}}_{\mathfrak{X}}(G/N)
\end{equation}
характеризует те пары $(G,N)$, где $N\trianglelefteqslant G$, для которых нахождение ${\mathfrak{X}}$-схемы группы~$G$ может быть заменено аналогичной задачей для меньшей группы~$G/N$. Подразумеваемое равенством~(\ref{1}) существование {\it некоторого} взаимно-од\-но\-знач\-ного соответствия между классами сопряженных ${\mathfrak{X}}$-максимальных подгрупп в группах~$G$ и~$G/N$ в~действительности означает, что
\begin{equation}\label{2}
HN/N\in\m_{\mathfrak{X}}(G/N)\text{ для любой }H\in\m_{\mathfrak{X}}(G),
\end{equation}
и отображение $H\mapsto HN/N$ индуцирует тем самым {\it естественное} вза\-им\-но-од\-но\-знач\-ное соответствие между этими классами.
В самом деле, пусть чертой обозначен канонический эпиморфизм $\overline{\phantom{G}}:G\rightarrow G/N.$ Хорошо известно (см., напр.,~\cite[лемма~2]{GR1}), что всякая ${\mathfrak{X}}$-подгруппа из $\overline{G}$ есть образ ${\mathfrak{X}}$-подгруппы из $G$. Отсюда $$\m_{\mathfrak{X}}(\overline{G})\subseteq\{\overline{H}\mid H\in \m_{\mathfrak{X}}(G)\}\quad\text{и}\quad{\mathrm{k}}_{\mathfrak{X}}(\overline{G})\leqslant {\mathrm{k}}_{\mathfrak{X}}(G).$$ Тем самым  из равенства~(\ref{1})
 следует равенство $$\m_{\mathfrak{X}}(\overline{G})=\{\overline{H}\mid H\in \m_{\mathfrak{X}}(G)\},$$ доказывающее справедливость утверждения (\ref{2}) и существование естественной биекции между классами сопряженных ${\mathfrak{X}}$-максимальных подгрупп в $G$ и~$G/N$.

Будем говорить, что
 \begin{itemize}
   \item {\it редукционная ${\mathfrak{X}}$-теорема верна  для пары $(G,N)$}, где $N\trianglelefteqslant G$, если выполнено~(\ref{1});
   \item {\it редукционная ${\mathfrak{X}}$-теорема верна для группы} $A$, если редукционная ${\mathfrak{X}}$-тео\-ре\-ма верна для любой пары $(G,N)$ такой, что $N\cong A$.
 \end{itemize}

Полагая $G=A$, видим, что редукционная ${\mathfrak{X}}$-теорема для группы $A$ влечет сопряженность ${\mathfrak{X}}$-максимальных подгрупп группы $A$:
$${\mathrm{k}}_{\mathfrak{X}}(A)={\mathrm{k}}_{\mathfrak{X}}(A/A)=1.$$
Виланд \cite[15.4]{Wie4} заметил, что, в свою очередь, редукционная ${\mathfrak{X}}$-теорема для $A$ вытекает предположения о сопряженности ${\mathfrak{X}}$-субмаксимальных подгрупп в $A$, и высказал гипотезу \cite[открытый вопрос к~15.4]{Wie4} о том, что сопряженность ${\mathfrak{X}}$-максимальных подгрупп  влечет сопряженность ${\mathfrak{X}}$-субмаксимальных подгрупп. В серии работ \cite{GR1,GR2,GR_surv,GR_SubmaxMinNonSolv,GRV,GRV_Reducktionssatz} авторам удалось подтвердить эту гипотезу~\cite[теорема~1]{GRV_Reducktionssatz}. Тем самым доказано, что
\begin{itemize}
  \item {\it редукционная ${\mathfrak{X}}$-теорема для группы $A$ равносильна тому, что ${\mathrm{k}}_{\mathfrak{X}}(A)=1$.}
\end{itemize}
Далее,
\begin{itemize}
  \item {\it условие  ${\mathrm{k}}_{\mathfrak{X}}(A)=1$ эквивалентно тому, что ${\mathrm{k}}_{\mathfrak{X}}(S)=1$ для любого композиционного фактора $S$ группы $A$,}
\end{itemize}
и для любой простой группы $S$ в терминах ее естественных арифметических параметров найдены необходимые и достаточные условия для того, чтобы выполнялось равенство ${\mathrm{k}}_{\mathfrak{X}}(S)=1$~\cite[теорема~1, приложение~A]{GRV_Reducktionssatz}. Например, {\it
для группы $S=\PSL_n(q)$, где $q$~--- порядок конечного поля характеристики $p\in\pi({\mathfrak{X}})$,  равенство ${\mathrm{k}}_{\mathfrak{X}}(S)=1$ равносильно тому, что \begin{itemize}
                                                           \item либо $S\in{\mathfrak{X}}$,
                                                           \item либо для любого $s\in\pi({\mathfrak{X}})$ если $s$ делит $|S|$, то $s$ делит~$q(q-1)$ и $s>n$.
                                                         \end{itemize}}
Таким образом, результаты работы~\cite{GRV_Reducktionssatz} можно интерпретировать, как описание всех пар $(G,N)$, для которых равенство~(\ref{1}) определяется только типом изоморфизма группы $N$.

Казалось бы, трудно ожидать, что редукционная ${\mathfrak{X}}$-теорема для пары $(G,N)$ всегда будет определяться одной только группой~$N$. Действительно, для группы $G$, обладающей нормальной подгруппой $N$, тип изоморфизма группы $N$ не определяет однозначно число ${\mathrm{k}}_{\mathfrak{X}}(G/N)$. Скажем, группа $$G=\PSL_2(7)\times \PGL_2(7)$$ обладает двумя нормальными подгруппами $N_1$ и $N_2$ такими, что $$N_1\cong N_2\cong\PSL_2(7),\quad\text{причем}\quad
G/N_1\cong\mathbb{Z}_2\times \PSL_2(7),\quad \text{а}\quad  G/N_2\cong\PGL_2(7).$$
Но для класса ${\mathfrak{X}}=\mathfrak{S}$ разрешимых групп легко показать (см., напр.,~\cite{Atlas}), что $${\mathrm{k}}_{\mathfrak{X}}(G/N_1)=3,\quad \text{а}\quad{\mathrm{k}}_{\mathfrak{X}}(G/N_2)=4.$$

На фоне данного примера теорема~\ref{Reduktionssatz_main} выглядит несколько неожиданно. Оказывается, то обстоятельство, что для данных группы~$G$ и ее нормальной подгруппы~$N$ мощность ${\mathfrak{X}}$-схемы не меняется при переходе от $G$ к $G/N$,  является внутренним свойством группы~$N$. Оно не зависит не только от особенностей вложения $N$ в~$G$, но и от самой $G$, поскольку влечет справедливость редукционной ${\mathfrak{X}}$-теоремы для~$N$. Именно этот факт является новым в теореме~\ref{Reduktionssatz_main} по сравнению с~\cite[теорема~1]{GRV_Reducktionssatz}. При этом если~\cite[теорема~1]{GRV_Reducktionssatz} была доказана путем сведения общей ситуации частному случаю ${\mathfrak{X}}=\mathfrak{G}_\pi$, изученному ранее в~\cite{MR,DpiCl,R2,R4,R5,Hall3',NumbCl} (см. также обзоры~\cite{Surv,GR_surv} и монографию~\cite{GuoBook1}), то утверждение ``только если'' в~теореме~\ref{Reduktionssatz_main} является новым даже для этого случая.

Наконец, учитывая, что в~\cite[теорема~1]{GRV_Reducktionssatz} получено описание всех {\it групп}, для которых верна редукционная ${\mathfrak{X}}$-теорема, получаем описание всех {\it пар} которых верна редукционная ${\mathfrak{X}}$-теорема:

\begin{corollary} \label{Description} Пусть ${\mathfrak{X}}$~--- полный класс групп. Редукционная ${\mathfrak{X}}$-теорема верна для пары $(G,N)$ тогда и только тогда, когда для любого композиционного фактора $S$ группы $N$ либо $S\in{\mathfrak{X}}$, либо пара $(S,{\mathfrak{X}})$ удовлетворяет одному из условий~{\rm I--VII} в {\rm \cite[приложение~A]{GRV_Reducktionssatz}.}
\end{corollary}

\subsection{Некоторые следствия}\label{Corollaries}

Поскольку ${\mathrm{k}}_{\mathfrak{X}}(G/N)\leqslant {\mathrm{k}}_{\mathfrak{X}}(G)$ для любой нормальной подгруппы $N$ группы~$G$, непосредственно из теоремы~\ref{Reduktionssatz_main} вытекает

\begin{corollary} \label{Inequality} Пусть ${\mathfrak{X}}$~--- полный класс групп и $N$~--- нормальная подгруппа группы конечной~$G$ такая, что ${\mathrm{k}}_{\mathfrak{X}}(N)>1$. Тогда ${\mathrm{k}}_{\mathfrak{X}}(G)>{\mathrm{k}}_{\mathfrak{X}}(G/N)$.
\end{corollary}

При этом, как показывает пример группы ${G=\PSL_2(7)\times \PGL_2(7)}$, рассмотренный выше, точное значение ${\mathrm{k}}_{\mathfrak{X}}(G/N)$ не только не определяется самими числами ${\mathrm{k}}_{\mathfrak{X}}(G)$ и ${\mathrm{k}}_{\mathfrak{X}}(N)$, но даже типом изоморфизма групп $G$~и~$N$.

Другим следствием теоремы~\ref{Reduktionssatz_main} оказывается существование в любой группе $G$ наибольшей нормальной подгруппы $R$ с тем свойством, что ${\mathrm{k}}_{\mathfrak{X}}(G)={\mathrm{k}}_{\mathfrak{X}}(G/R)$.

\begin{corollary} \label{Radical} Пусть ${\mathfrak{X}}$~--- полный класс. Для произвольной конечной группы $G$ рассмотрим подгруппу
$$R=\langle N\mid N\trianglelefteqslant G\quad \text{и}\quad{\mathrm{k}}_{\mathfrak{X}}(G)={\mathrm{k}}_{\mathfrak{X}}(G/N)\rangle.$$ Она обладает следующими свойствами:
\begin{itemize}
  \item[$(i)$] ${\mathrm{k}}_{\mathfrak{X}}(G)={\mathrm{k}}_{\mathfrak{X}}(G/R);$
   \item[$(ii)$] если $N\trianglelefteqslant G$ и $N\leqslant R$, то ${\mathrm{k}}_{\mathfrak{X}}(G)={\mathrm{k}}_{\mathfrak{X}}(G/N);$
  \item[$(iii)$] если $\overline{G}=G/R$, то ${\mathrm{k}}_{\mathfrak{X}}(\overline{G})={\mathrm{k}}_{\mathfrak{X}}(\overline{G}/\overline{N})$ влечет $\overline{N}=1$ для любой $\overline{N}\trianglelefteqslant \overline{G}$.
\end{itemize}
\end{corollary}
Ввиду теоремы~\ref{Reduktionssatz_main}, подгруппа $R\leqslant G$, о которой идет речь в следствии~\ref{Radical}, совпадает с $\D_{\mathfrak{X}}$-радикалом  группы $G$, где, как и в \cite{GR1,GR2,GRV_Reducktionssatz}, через $\D_{\mathfrak{X}}$ обозначен класс конечных групп в которых  все ${\mathfrak{X}}$-максимальные подгруппы сопряжены. Класс $\D_{\mathfrak{X}}$ замкнут относительно взятия нормальных подгрупп гомоморфных образов и расширений\footnote{См.~\cite[следствие~1]{GRV_Reducktionssatz}. Отметим также, что в силу неравенства ${\mathrm{k}}_{\mathfrak{X}}(G/N)\leqslant {\mathrm{k}}_{\mathfrak{X}}(G)$ данное утверждение является частным случаем теоремы~\ref{Reduktionssatz_main}.} и, в частности, является классом Фиттинга (см. определение в \cite{DH,GuoBook,GuoBook1}), и поэтому любая группа обладает $\D_{\mathfrak{X}}$-радикалом. При этом, вообще говоря, $\D_{\mathfrak{X}}$ не является полным классом, поскольку может оказаться не замкнутым относительно взятия подгрупп (см.~\cite[теорема~1.7]{VMR}).

Факторгруппу $G/R$ назовем {\it полной редукцией над ${\mathfrak{X}}$ группы $G$}, а саму подгруппу $R=\langle N\mid N\trianglelefteqslant G\quad \text{и}\quad{\mathrm{k}}_{\mathfrak{X}}(G)={\mathrm{k}}_{\mathfrak{X}}(G/N)\rangle$~--- {\it ядром этой редукции}.
Группу $G$ назовем {\it вполне редуцированной над ${\mathfrak{X}}$}, если ядро ее полной ${\mathfrak{X}}$-редукции равно единице. Результаты данной работы сводят задачу о нахождении ${\mathfrak{X}}$-схемы к случаю вполне редуцированных групп.

Для группы $G$ обозначим через $\mathrm{om}_{\mathfrak{X}}(G)$ множество надгрупп ${\mathfrak{X}}$-максимальных подгрупп, т.\,е.
$$
\mathrm{om}_{\mathfrak{X}}(G)=\{K\leqslant G\mid \m_{\mathfrak{X}}(K)\cap\m_{\mathfrak{X}}(G)\ne\varnothing\}.
$$
Менее очевидным утверждением, вытекающим из теоремы~\ref{Reduktionssatz_main} и основного результата работы~\cite{VMR}, будет

\begin{corollary} \label{Overgroups} Пусть ${\mathfrak{X}}$~--- полный класс групп и $N$~--- нормальная подгруппа группы конечной~$G$. Следующие утверждения эквивалентны.
\begin{itemize}
  \item[$(i)$] ${\mathrm{k}}_{\mathfrak{X}}(G)={\mathrm{k}}_{\mathfrak{X}}(G/N);$
  \item[$(ii)$] ${\mathrm{k}}_{\mathfrak{X}}(K)={\mathrm{k}}_{\mathfrak{X}}\left(K/(K\cap N)\right)$ для всех $K\in \mathrm{om}_{\mathfrak{X}}(G)$.
\end{itemize}
\end{corollary}

\subsection{ Категория групп и ${\mathfrak{X}}$-изосхематизмов}\label{Iso}

Переведем теорему~\ref{Reduktionssatz_main} на язык гомоморфизмов.
%Тем более важным представляется то обстоятельство, что эпиморфизмы, сохраняющие полную информацию об ${\mathfrak{X}}$-схеме существуют и, более того, как  показано в данной работе, их можно полностью описать.
Назовем эпиморфизм групп $$\phi:G \rightarrow G^*$$  {\it изосхематизмом над ${\mathfrak{X}}$} или {\it ${\mathfrak{X}}$-изосхематизмом}, если он переводит  ${\mathfrak{X}}$-схему группы~$G$ (любую или некоторую) в ${\mathfrak{X}}$-схему группы $G^*$. Теорема~\ref{Reduktionssatz_main} эквивалентна следующему утверждению.

\begin{itemize}
  \item {\it Эпиморфизм $\phi$ является ${\mathfrak{X}}$-изосхематизмом тогда и только тогда, когда $${\mathrm{k}}_{\mathfrak{X}}(\ker\phi)=1.$$}
\end{itemize}

Из сказанного выше следует, что свойство эпиморфизма $\phi:G \rightarrow G^*$ быть или не быть ${\mathfrak{X}}$-изо\-схе\-ма\-тиз\-мом полностью определяется группами самими $G$ и~$G^*$ и не зависит от конкретного отображения. Сформулируем это явным образом.

\begin{propos} \label{IsoschematismDefinitions} Пусть $G$~--- конечная группа, а $G^*$~--- ее эпиморфный образ. Для полного класса ${\mathfrak{X}}$ следующие утверждения эквивалентны.
\begin{itemize}
  \item[$(i)$] Существует ${\mathfrak{X}}$-изосхематизм $\phi:G \rightarrow G^*$.
  \item[$(ii)$] Любой эпиморфизм $\phi:G \rightarrow G^*$ является  ${\mathfrak{X}}$-изосхематизмом.
  \item[$(iii)$] ${\mathrm{k}}_{\mathfrak{X}}(G)={\mathrm{k}}_{\mathfrak{X}}(G^*)$.
  \end{itemize}
\end{propos}

Отметим, что ядра любых двух ${\mathfrak{X}}$-изосхематизмов из $G$ на $G^*$ могут быть не изоморфными, хотя из теоремы Жордана--Гёльдера следует они имеют один и тот же набор композиционных факторов.

Тот факт, что существует ${\mathfrak{X}}$-изосхематизм из~$G$ на~$G^*$, будем записывать, как
$$G\underset{{\mathfrak{X}}}{\twoheadrightarrow} G^*.$$
Этот же символ будет использоваться в записи
$$\phi:G\underset{{\mathfrak{X}}}{\twoheadrightarrow} G^*,$$
означающей, что отображение $\phi$ является ${\mathfrak{X}}$-изосхематизмом из $G$ на~$G^*$.

Можно рассматривать $\underset{{\mathfrak{X}}}{\twoheadrightarrow}$, как отношение между группами. Оно очевидно рефлексивно и транзитивно, но не симметрично. Симметризуем его. Скажем, что две группы $G_1$ и $G_2$ {\it изосхемны над ${\mathfrak{X}}$} или {\it ${\mathfrak{X}}$-изосхемны}, и будем писать
$$G_1\underset{{\mathfrak{X}}}{\equiv} G_2,$$ если из $G_1$ и $G_2$ существуют ${\mathfrak{X}}$-изосхематизмы на одну и ту же группу:
$$\xymatrix{
  G_1 %\ar@{~}[rr]
  \ar@{->>}[dr]_{\mathfrak{X}}
                &
\underset{{\mathfrak{X}}}{\equiv} &    G_2 \ar@{->>}[dl]^{\mathfrak{X}}    \\
                & G_0                 }
$$
Отношение $\underset{{\mathfrak{X}}}{\equiv}$ очевидно рефлексивно и симметрично. В действительности, оно задает отношение эквивалентности на группах, а его транзитивность вытекает из теоремы~\ref{Reduktionssatz_main}. Это отношение позволяет описать категорию групп и ${\mathfrak{X}}$-изосхематизмов.

\begin{corollary} \label{IsoschematismsCategory} Соотношение $G_1\underset{{\mathfrak{X}}}{\equiv} G_2$ для конечных групп $G_1$ и $G_2$ равносильно тому, что полные редукции над ${\mathfrak{X}}$ этих групп изоморфны. Отношение $\underset{{\mathfrak{X}}}{\equiv}$ является отношением эквивалентности между конечными группами. Каждый класс эквивалентности содержит единственную с точностью до изоморфизма вполне редуцированную над~${\mathfrak{X}}$ группу, которая является универсально притягивающим объектом\footnote{См. определение в \cite[гл.~1, \S7]{Lang}. Отметим, что ${\mathfrak{X}}$-изосхематизмы в качестве морфизмов в данной категории рассматриваются с точностью до композиции с автоморфизмами групп.} в этом классе как подкатегории в категории всех конечных групп и ${\mathfrak{X}}$-изо\-схе\-ма\-тиз\-мов.
\end{corollary}

Следствие~\ref{Overgroups} на языке гомоморфизмов  можно переформулировать, как следующее утверждение.
 \begin{itemize}
   \item {\it Пусть $\phi$~--- ${\mathfrak{X}}$-изосхематизм, определенный на группе $G$ и $K$~--- над\-груп\-па  ${\mathfrak{X}}$-мак\-си\-маль\-ной подгруппы из~$G$. Тогда ограничение $\phi$ на $K$  является ${\mathfrak{X}}$-изо\-схе\-ма\-тиз\-мом $K \underset{{\mathfrak{X}}}{\twoheadrightarrow} K^\phi$.}

 \end{itemize}

\section{ Обозначения и предварительные леммы}

Используемые нами обозначения из теории групп стандартны и могут быть найдены в~\cite{Atlas,Bray,DH,GuoBook}. Для натурального числа $n$ через $\pi(n)$ обозначается множество его простых делителей, а для группы $G$ полагаем $\pi(G)=\pi(|G|)$.  Для фиксированного множества $\pi\subseteq\mathbb{P}$ простых чисел и полного класса ${\mathfrak{X}}$ конечных групп мы используем следующие менее стандартные обозначения.
\begin{itemize}
\item[] $\Omega/G$ для случая, когда группа $G$ действует на множестве $\Omega$ --- множество орбит этого действия.
\item[] $|\Omega:G|$ для случая, когда группа $G$ действует на множестве $\Omega$ --- число орбит этого действия, т.е. $|\Omega:G|=|\Omega/G|$.
\item[] $\Hall_\pi(G)$ --- множество $\pi$-холловых подгрупп группы $G$.
\item[] $\Hall_{\mathfrak{X}}(G)$ --- множество {\it ${\mathfrak{X}}$-холловых подгрупп} группы $G$, т.\,е. таких ${\mathfrak{X}}$ подгрупп, индекс которых не делится ни на какие числа из~$\pi({\mathfrak{X}})$.
\item[] $\m_{\mathfrak{X}}(G)$ --- множество {\it ${\mathfrak{X}}$-максимальных подгрупп} группы $G$.
\item[] ${\mathrm{k}}_{\mathfrak{X}}(G)$ --- число классов сопряженности ${\mathfrak{X}}$-максимальных подгрупп группы $G$, т.\,е. $${\mathrm{k}}_{\mathfrak{X}}(G)=|\m_{\mathfrak{X}}(G):G|$$ для действия группы $G$ сопряжениями на множестве $\m_{\mathfrak{X}}(G)$.
\item[] $\h_{\mathfrak{X}}(G)$ --- число классов сопряженности ${\mathfrak{X}}$-холловых подгрупп группы $G$, т.\,е. $$\h_{\mathfrak{X}}(G)=|\Hall_{\mathfrak{X}}(G):G|$$ для действия группы $G$ сопряжениями на множестве $\Hall_{\mathfrak{X}}(G)$.
\item[] $\E_{\mathfrak{X}}$ ---  класс всех конечных групп $G$ таких, что $\Hall_{\mathfrak{X}}(G)\ne\varnothing$ (эквивалентно, $\h_{\mathfrak{X}}(G)\ge 1$).
 \item[] $\C_{\mathfrak{X}}$ ---  класс всех конечных групп $G$ таких, что $\h_{\mathfrak{X}}(G)=1$.
 \item[] $\D_{\mathfrak{X}}$ ---  класс всех конечных групп $G$ таких, что ${\mathrm{k}}_{\mathfrak{X}}(G)=1$.
 \item[] $\M_{\mathfrak{X}}$ ---  класс всех конечных групп $G$ таких, что $\m_{\mathfrak{X}}(G)=\Hall_{\mathfrak{X}}(G)$ (эквивалентно, ${\mathrm{k}}_{\mathfrak{X}}(G)=\h_{\mathfrak{X}}(G)$.
\end{itemize}

Обозначения $\E_{\mathfrak{X}}$, $\C_{\mathfrak{X}}$ и $\D_{\mathfrak{X}}$ обобщают обозначения $\E_\pi$, $\C_\pi$ и $\D_\pi$, введенные Ф.\,Холлом \cite{Hall}, и совпадают с ними в случае, когда ${\mathfrak{X}}=\mathfrak{G}_\pi$~--- класс всех $\pi$-групп. Из теоремы Силова следует, что $\D_{\mathfrak{X}}=\C_{\mathfrak{X}}\cap \M_{\mathfrak{X}}$ и включения между классами $\E_{\mathfrak{X}}$, $\C_{\mathfrak{X}}$, $\M_{\mathfrak{X}}$ и $\D_{\mathfrak{X}}$ отражены на следующей диаграмме:
$$
\xymatrix{
         & \E_{\mathfrak{X}} \ar@{-}[dr]^{} & &\\
  \C_{\mathfrak{X}} \ar@{-}[ur]_{} %\ar@{-}[rr]^{}
  \ar@{-}[dr]_{}
                &  &   \M_{\mathfrak{X}} \ar@{-}[dl]^{}    \\
                &                \D_{\mathfrak{X}}=\C_{\mathfrak{X}}\cap \M_{\mathfrak{X}}                 }
%\xymatrix{
%                &  \ar[dr]^{}             \\
%  \ar[ur]^{} \ar[rr]^{} & &             }
 $$
 В случае ${\mathfrak{X}}=\mathfrak{G}_\pi$, подобно тому, как, мы пишем $\E_\pi$ вместо $\E_{\mathfrak{X}}$, мы будем использовать естественные обозначения ${\mathrm{k}}_\pi(G)$ и $\h_\pi(G)$ для числа классов сопряженности максимальных $\pi$-подгрупп и $\pi$-холловых подгрупп группы~$G$.  Будем говорить, что $n\in\mathbb{N}$~--- {\it $\pi$-число}, если все его простые делители принадлежат~$\pi$.

 \begin{lemma}\label{FrattiniSubgr} {\em \cite[лемма~2]{GR1}} Пусть ${\mathfrak{X}}$~--- полный класс. Допустим, $$\phi:G\rightarrow G_0$$
  является гомоморфизмом групп и предположим, что $K\in{\mathfrak{X}}$ для некоторой подгруппы $K\leq G^\phi$. Тогда $K=H^\phi$ для некоторой ${\mathfrak{X}}$-подгруппы $H\leq G$. В частности, $\m_{\mathfrak{X}}(G^\phi)\subseteq\m_{\mathfrak{X}}(G)^\phi$.
\end{lemma}

Через ${\mathfrak{X}}'$ обозначим класс всех групп $G$ таких, что $\m_{\mathfrak{X}}(G)=\{1\}$. Группу называем {\it ${\mathfrak{X}}$-отделимой}, если она обладает (суб)нормальным рядом, каждый фактор которого~--- ${\mathfrak{X}}$- или ${\mathfrak{X}}'$-группа.

В следующей лемме собраны некоторое известные результаты о поведении ${\mathfrak{X}}$-максимальных и ${\mathfrak{X}}$-холловых подгрупп.

 \begin{lemma}\label{HallSubgroup}\label{Lifting}\label{Bijection} Пусть $N$~--- нормальная подгруппа группы $G$. Тогда справедливы следующие утверждения.
 \begin{itemize}
   \item[$(i)$] {\em \cite[лемма~1]{Hall}.}
Если  $H\in\Hall_{\mathfrak{X}}(G)$, то $${H\cap N\in \Hall_{\mathfrak{X}}(N)}\quad\text{и}\quad HN/N\in \Hall_{\mathfrak{X}}(G/N).$$
   \item[$(ii)$] {\em \cite[лемма~2.1(e)]{NumbCl}.}
 Допустим, $G/N\in{\mathfrak{X}}$. Тогда для $H\in\Hall_{\mathfrak{X}}(N)$ в том и только в том случае найдется $K\in\Hall_{\mathfrak{X}}(G)$, для которой с $K\cap N$, когда $H^S=H^G$ $($другими словами, когда класс $H^S\in\Hall_{\mathfrak{X}}(N)/N$ инвариантен относительно действия группы $G$ на множестве $\Hall_{\mathfrak{X}}(N)/N$ сопряжениями$)$.
   \item[$(iii)$] {\rm \cite[12.9]{Wie4}.} {\it Допустим, $N$~--- ${\mathfrak{X}}$-отделимая группа Тогда  ${\mathrm{k}}_{\mathfrak{X}}(G)={\mathrm{k}}_{\mathfrak{X}}(G/N)$. В~частности,  $G\in \D_{\mathfrak{X}}$ тогда и только тогда, когда $G/N\in \D_{\mathfrak{X}}$.}
 \end{itemize}

\end{lemma}

\begin{lemma}\label{RedSatz} {\em ~\cite[теорема~1]{GRV_Reducktionssatz}.} Пусть $N$~--- нормальная подгруппа группы~$G$ и ${{\mathrm{k}}_{\mathfrak{X}}(N)=1}$. Тогда ${\mathrm{k}}_{\mathfrak{X}}(G)={\mathrm{k}}_{\mathfrak{X}}(G/N)$.
\end{lemma}

 Если $S$ и $H$~--- подгруппы группы $G$, обозначим через $\Aut_H(S)$ {\it группу $H$-ин\-ду\-ци\-ро\-ван\-ных автоморфизмов группы $S$}, т.\,е. образ в $\Aut(S)$ гомоморфизма $$\alpha_H:N_H(S)\rightarrow \Aut(S),$$ который сопоставляет любому элементу $x\in N_H(S)$ автоморфизм группы $S$, заданный правилом $s\mapsto s^x=x^{-1}sx$. Ядро этого гомоморфизма равно $C_H(S)$, поэтому имеем
 $$
 \Aut_H(S)\cong N_H(S)/C_H(S).
 $$
 Если при этом $H\le K\le G$, то определенный выше гомоморфизм $\alpha_H$ является ограничением на $N_H(S)$ соответствующего гомоморфизма $\alpha_K:N_K(S)\rightarrow\Aut_K(S)$. Следовательно, $\Aut_H(S)\le \Aut_K(S)$. %ввиду того, что
% $$
% \Aut_H(S)\cong N_H(S)/C_H(S)\cong (H\cap N_K(S))C_K(S)/C_K(S)\le  N_K(S)/C_K(S)\cong \Aut_K(S).
%$$

\begin{lemma}\label{Constrain} Пусть $S$ --- простая неабелева субнормальная подгруппа группы $G$ и $H\in\Hall_{\mathfrak{X}}(S)$, причем группа $\Aut_G(S)$ стабилизирует класс сопряженности подгруппы $H$ в $S$ $($т.\,е.~$H^{\Aut_G(S)}=H^S$$).$ Возьмем произвольную систему представителей $g_1,\dots,g_n$
правых смежных классов группы $G$ по $N_G(S)$. Пусть $$M=\langle S^{g_i}\mid
i=1,\dots, n\rangle\text{ и }V= \langle H^{g_i}\mid i=1,\dots, n\rangle.$$ Тогда
$M=\langle S^G \rangle$~--- минимальная нормальная подгруппа группы $G$, и имеют место следующие утверждения:
\begin{itemize}
  \item[$(i)$] $V^G=V^M$;
  \item[$(ii)$]$V\in\Hall_{\mathfrak{X}}(M)$;.
  \item[$(iii)$] если $G/M\in {\mathfrak{X}}$, то $V=K\cap M$ и $H=K\cap S$ для некоторой $K\in\Hall_{\mathfrak{X}}(G)$.
\end{itemize}
\end{lemma}
\noindent{\it Доказательство.}
 Ввиду выбора элементов $g_1,\dots,g_n$ для любого $g\in G$ имеем
$$S^g\in\{S^{g_i}\mid i=1,\dots, n\}$$ и поэтому $M=\langle S^G \rangle\trianglelefteqslant G$.
Заметим также, что ввиду простоты и субнормальности подгруппы $S$ имеем
$[S^{g_i},S^{g_j}]=1$ при $i\ne j$.

Пусть $g\in G$. Найдутся подстановка
$\sigma\in\Sym_n$ и элементы $x_1,\dots,x_n\in N_G(S)$ такие, что
$g_ig=x_ig_{i\sigma}$. Рассмотрим автоморфизмы $\gamma_i\in\Aut_G(S)$, задаваемые  правилом $\gamma_i:s\mapsto s^{x_i}$. По условию
$H^{x_i}=H^{\gamma_i}=H^{s_i}$ для некоторого $s_i\in S$. Положим
$a_i=s_{i\sigma^{-1}}^{g_i}$ и $a=a_1\dots a_n$. Ясно, что $a\in M$.  Покажем, что $V^g=V^a$ и тем самым установим равенство $V^G=V^M$.

Из определения следует, что
$a_i\in S^{g_i}$ и $H^{g_ia}=H^{g_ia_i}$. Имеем
\begin{multline*}
V^g= \langle H^{g_ig}\mid i=1,\dots, n\rangle= \langle H^{x_ig_{i\sigma}}\mid
i=1,\dots, n\rangle=\\ \langle H^{s_ig_{i\sigma}}\mid i=1,\dots, n\rangle=
 \langle H^{s_{i\sigma^{-1}}g_{i}}\mid i=1,\dots, n\rangle=\\ \langle
H^{g_is_{i\sigma^{-1}}^{g_i}}\mid i=1,\dots, n\rangle=\langle H^{g_ia_i}\mid
i=1,\dots, n\rangle=\\
 \langle H^{g_ia}\mid i=1,\dots, n\rangle=V^a.
\end{multline*}
Утверждение $(i)$ доказано. Далее,
$V$ является прямым произведением ${\mathfrak{X}}$-групп $H^{g_i}$, $i=1,\dots,n$, поэтому~$V\in{\mathfrak{X}}$. Кроме того, так как ${H\in\Hall_{\mathfrak{X}}(S)}$, число
$$
|M:V|=\prod_{i=1}^n|S^{g_i}:H^{g_i}|=|S:H|^n
$$
не делится ни на какие числа из $\pi({\mathfrak{X}})$. Поэтому ${V\in\Hall_{\mathfrak{X}}(M)}$, что доказывает~$(ii)$.  Наконец, из $(i)$ и леммы~\ref{Lifting} получаем~$(iii)$.
\qed\par\medskip

\begin{lemma}\label{AutK=AutNK}
 {\it Пусть нормальная подгруппа $N$ группы $G$ является прямым произведением неабелевых простых групп и $S$~--- одна из них. Предположим, что $G=KN$ для некоторой подгруппы $K$. Тогда

 \begin{itemize}
   \item[$(i)$] $N_G(S)=NN_K(S)$
   \item[$(ii)$] $\Aut_G(S)=\Inn(S)\Aut_K(S)$
   %\item[$(iii)$]
 \end{itemize}
 }
\end{lemma}

\noindent{\it Доказательство.}
Пусть $N=S_1\times S_2\times\dots\times S_n$ и $S_1=S$. Тогда $$N\le N_G(S)\text{ и }S_2\times\dots\times S_n=C_N(S)\le C_G(S).$$ Поэтому
$$
N_G(S)=NN_K(S),
$$
как и утверждается в $(i)$. Пусть $$\alpha:N_G(S)\rightarrow \Aut(S)$$ обозначает естественный гомоморфизм, индуцированный сопряжениями. Его ядро равно $C_G(S)$.   Имеем $S^\alpha=\Aut_S(S)=\Inn(S)$, %откуда $$\Inn(S)\Aut_K(S)\le \Aut_G(S).$$ Наконец,
$N=SC_N(S)$, поэтому $N^\alpha= \Inn(S)$ и
$$
%\begin{multline*}
\Aut_G(S)= N_G(S)^\alpha= N^\alpha N_K(S)^\alpha = %=|NC_G(S)/C_G(S)||N_K(S)C_G(S)/C_G(S)|=\\|N/C_N(S)||N_K(S)/C_K(S)|=
\Inn(S)\Aut_K(S).
%\end{multline*}
$$
Утверждение~$(ii)$ также доказано.
\qed\par\medskip

\medskip

Ключевую роль в доказательстве теоремы~\ref{Reduktionssatz_main} играет теорема о числе классов сопряженных $\pi$-холловых подгрупп в простых группах, доказанная в~\cite{NumbCl}. Нам она понадобится в следующем уточненном виде.

\begin{lemma}\label{ClassNumberpi}{\rm \cite[теорема~1.1]{NumbCl} } {\it
Пусть  $S$~--- простая конечная группа, обладающая холловой $\pi$-подгруппой для некоторого множества $\pi$ простых чисел. Тогда имеет место одно из следующих утверждений.
\begin{itemize}
  \item[$(i)$]   $2\notin\pi$ и ${\rm h}_\pi(S)=1$.
  \item[$(ii)$]   $3\notin\pi$ и ${\rm h}_\pi(S)\in\{1,2\}$.
  \item[$(iii)$]   $2,3\in\pi$ и ${\rm h}_\pi(S)\in\{1,2,3,4,9\}$.
  \end{itemize}
}
\end{lemma}

\begin{lemma}\label{ClassNumber=9} {\rm \cite[Лемма 12]{GR2}} {\it Пусть ${\mathfrak{X}}$~--- полный класс. Положим $\pi=\pi({\mathfrak{X}})$. Предположим также, что  ${\rm h}_\pi(S)=9$. Тогда ${\rm h}_{\mathfrak{X}}(S)$ совпадает с одним из чисел $0$, $1$ или~$9$.}\end{lemma}

Из лемм~\ref{ClassNumberpi} и~\ref{ClassNumber=9} вытекает

\begin{lemma}\label{ClassNumberX} {\it
Пусть  $S$~--- простая конечная группа. Тогда имеет место одно из следующих утверждений.
\begin{itemize}
  \item[$(i)$]   $2\notin\pi({\mathfrak{X}})$ и ${\rm h}_{\mathfrak{X}}(S)\in\{0,1\}$.
  \item[$(ii)$]   $3\notin\pi({\mathfrak{X}})$ и ${\rm h}_{\mathfrak{X}}(S)\in\{0,1,2\}$.
  \item[$(iii)$]   $2,3\in\pi({\mathfrak{X}})$ и ${\rm h}_{\mathfrak{X}}(S)\in\{0,1,2,3,4, 9\}$.%
  \end{itemize}
}
\end{lemma}

Предположим, что для простой группы $S$ выполнено равенство  ${{\rm h}_{\mathfrak{X}}(S)=9}$. Поскольку ${{\rm h}_{\mathfrak{X}}(S)\leqslant{\rm h}_\pi(S)}$ для ${\pi=\pi({\mathfrak{X}})}$, по лемме~\ref{ClassNumberpi} получаем ${{\rm h}_{\mathfrak{X}}(S)={\rm h}_\pi(S)}$ и следовательно $\Hall_{\mathfrak{X}}(S)=\Hall_\pi(G)$. Теперь, используя информацию из \cite[леммы~2.3, 3.1, 4.4, и 8.1]{NumbCl}, мы получаем следующую лемму о строении ${\mathfrak{X}}$-хол\-ло\-вых подгрупп группы~$S$.

\begin{lemma}\label{ClassNumberX=9}
Пусть  $S$~--- простая конечная группа и  ${\rm h}_{\mathfrak{X}}(S)=9$ для некоторого полного класса~${\mathfrak{X}}$. Тогда справедливы следующие утверждения.
\begin{itemize}
  \item[$(i)$] $S\cong \PSp_{2n}(q)\cong \PSp(V)$, где $q$~--- степень простого числа $p\notin \pi({\mathfrak{X}})$, а $V$~--- ассоциированное с $\PSp_{2n}(q)$ векторное  пространство размерности~$2n$ над полем~$\mathbb{F}_q$ с невырожденной кососимметрической формой.
  \item[$(ii)$] $\pi({\mathfrak{X}})\cap\pi(S)=\pi(q^2-1)$ и
  \begin{itemize}
    \item[$\bullet$] либо $\pi({\mathfrak{X}})\cap\pi(S)=\{2,3\}$ и $n\in\{5,7\}$,
    \item[$\bullet$] либо $\pi({\mathfrak{X}})\cap\pi(S)=\{2,3,5\}$ и $n=7.$
  \end{itemize}
  \item[$(iii)$] Любая холлова $\pi({\mathfrak{X}})$-подгруппа группы $\PSp_{2n}(q)$ содержится в стабилизаторе~$M$ разложения пространства $V$ с в ортогональную сумму
       $$V=V_1\perp \dots\perp V_n$$ невырожденных изометричных подпространств размерности~$2$. Существует подгруппа $A\trianglelefteqslant M$, такая, что $A=L_1\dots L_n$, где $L_i\cong \Sp(V_i)\cong \Sp_2(q)\cong SL_2(q)$, $[L_i,L_j]=1$, $i,j=1,\dots,n$, $i\ne j$, и $M/A\cong \Sym_n$.
  \item[$(iv)$] ${\rm h}_{\mathfrak{X}}(\Sym_n)=1$. При этом
  \begin{itemize}
    \item[$\bullet$] если $\pi({\mathfrak{X}})\cap\pi(S)=\{2,3\}$, то ${\mathfrak{X}}$-холлова подгруппа группы $\Sym_n$ изоморфна $\Sym_4$ для $n=5$ и $\Sym_3\times \Sym_4$ для $n=7$;
    \item[$\bullet$] если $\pi({\mathfrak{X}})\cap\pi(S)=\{2,3,5\}$, то ${\mathfrak{X}}$-холлова подгруппа группы $\Sym_n=\Sym_7$ изоморфна $\Sym_6$ $($в частности, $\Sym_m\in{\mathfrak{X}}$ при $m\leqslant 6$$)$.
  \end{itemize}
  \item[$(v)$] ${\rm h}_{\mathfrak{X}}(\Sp_2(q))=3$. При этом
  \begin{itemize}
    \item[$\bullet$] если $\pi({\mathfrak{X}})\cap\pi(S)=\{2,3\}$, то все ${\mathfrak{X}}$-хол\-ло\-вы подгруппы в $\Sp_2(q)\cong \SL_2(q)$ разрешимы, $\Sp_2(q)$ содержит один класс сопряженных ${\mathfrak{X}}$-хол\-ло\-вых подгрупп, изоморфных обобщенной группе кватернионов порядка $48$, и два класса ${\mathfrak{X}}$-хол\-ло\-вых подгрупп, изоморфных $2.\Sym_4$;
    \item[$\bullet$] если $\pi({\mathfrak{X}})\cap\pi(S)=\{2,3,5\}$, то  группа $\Sp_2(q)\cong \SL_2(q)$ содержит один сопряженности класс разрешимых ${\mathfrak{X}}$-хол\-ло\-вых подгрупп, изоморфных обобщенной группе кватернионов порядка~$120$, и два класса ${\mathfrak{X}}$-хол\-ло\-вых подгрупп, изоморфных~$\SL_2(5)\cong 2.\Alt_5$.
  \end{itemize}
  \item[$(vi)$] Число неподвижных точек любой подгруппы $G\leq\Aut(S)$ при ее действии  на множестве $\Hall_{\mathfrak{X}}(S)/S$ равно $1$ или~$9$.
\end{itemize}
\end{lemma}

\section{Аргумент Фраттини для ${\mathfrak{X}}$-холловых подгрупп}

Основная цель данного раздела в том, чтобы доказать следующее утверждение:

\begin{propos}\label{XFrattini} {\em }
Пусть $G$ группа обладает нормальной подгруппой $A$ такой, что $A=KN$ для некоторой нормальной в~$G$ подгруппы $N$, являющейся прямым про\-из\-ве\-де\-нием неабелевых простых групп, и некоторой $K\in\Hall_{\mathfrak{X}}(G)$. Тогда существует $L\in\Hall_{\mathfrak{X}}(A)$ такая, что $G=AN_G(L)$.
\end{propos}

 \noindent{\it Доказательство.}%[предложения~\ref{XFrattini}]
 Пусть $\pi=\pi({\mathfrak{X}})$. Так как $A$ содержит $H\in\Hall_{\mathfrak{X}}(G)$, индекс $|{G:A}|$ является $\pi'$-числом. Так как подгруппа $A$ нормальна в $G$, имеем $\Hall_{\mathfrak{X}}(G)=\Hall_{\mathfrak{X}}(A)$ и $G/A$~--- $\pi'$-группа.

 Пусть
 $$N=S_1\times\dots\times S_n,$$ где $S_1,\dots,S_n$~--- неабелевы простые группы. Нам потребуется установить некоторые факты о холловых ${\mathfrak{X}}$-подгруппах групп $S_1,\dots, S_n$, о классах сопряженности таких подгрупп и о действии на этих классах групп $G$-индуцированных автоморфизмов. Для этого зафиксируем произвольно $S\in\{S_1,\dots,S_n\}$. Так как $K\cap S\in\Hall_{\mathfrak{X}}(S)$, имеем $S\in\E_{\mathfrak{X}}$ и по лемме~\ref{ClassNumberX} $$\h_{\mathfrak{X}}(S)\in\{1,2,3,4,9\}.$$

 Пусть $\Omega$~--- множество неподвижных точек группы $\Aut_A(S)$ на множестве $\Hall_{\mathfrak{X}}(S)/S$ классов сопряженности ${\mathfrak{X}}$-холловых подгрупп группы $S$, т.\,е.
  \begin{multline*}
 \Omega=\{H^S\mid H\in\Hall_{\mathfrak{X}}(S)\text{ и }\\
  \text{для любого }a\in A\text{ существует }x\in S\text{ такой, что }H^a=H^x\}.
 \end{multline*}
  Заметим, что $\Omega\ne\varnothing$, поскольку $(K\cap S)^S\in\Omega$. В самом деле, $N\le N_A(S)$, поэтому $$N_A(S)=N_{KN}(S)=N_K(S)N.$$ При этом класс сопряженности $(K\cap S)^S$ инвариантен относительно обеих групп $N_K(S)$ и $N$. Поэтому он инвариантен относительно $N_A(S)$ и относительно $\Aut_A(S)$ и, значит, принадлежит~$\Omega$.

  Так как $A\trianglelefteqslant G$, имеем  $N_A(S)\trianglelefteqslant N_G(S)$. Отсюда $\Aut_A(S)\trianglelefteqslant \Aut_G(S)$ и, следовательно, группа $\Aut_G(S)$ действует на $\Omega$. Мы утверждаем, что
 \begin{itemize}
   \item[$1^\circ$] {\it длина некоторой орбиты группы $\Aut_G(S)$ на $\Omega$ является $\pi$-числом.}
 \end{itemize}

 В самом деле, $|\Omega|\le \h_{\mathfrak{X}}(S)$, а длина любой орбиты группы $\Aut_G(S)$ на $\Omega$ не превосходит $|\Omega|$. Если $2\notin\pi$ или $3\notin\pi$, то из леммы~\ref{ClassNumberX} следует, что длина любой орбиты группы $\Aut_G(S)$ на $\Omega$ является $\pi$-числом. Поэтому считаем, что $2,3\in\pi$. Теперь если $\h_{\mathfrak{X}}(S)\le 4$, то снова,  длина любой орбиты группы $\Aut_G(S)$ на $\Omega$ является $\pi$-числом. Поэтому считаем, что $\h_{\mathfrak{X}}(S)=9$. Из леммы~\ref{ClassNumberX=9} следует, что $$|\Omega|\in\{1,9\}.$$ Случай $|\Omega|=1$ очевиден. Из того, что единственными не $\{2,3\}$-числами, не пре\-вос\-хо\-дя\-щи\-ми~$9$, являются $5$ и~$7$, легко выводится также, что в любом разбиении числа $9$ в сумму натуральных слагаемых присутствует слагаемое, являющееся $\{2,3\}$-числом и, следовательно, $\pi$-числом. Поэтому среди орбит, на которые распадается $\Omega$ относительно действия группы $\Aut_G(S)$, есть орбита, длина которой является $\pi$-числом.

 Утверждение $1^\circ$ можно усилить следующим образом:
 \begin{itemize}
   \item[$2^\circ$] {\it если длина некоторой орбиты группы $\Aut_G(S)$ на $\Omega$ является $\pi$-числом, то эта длина равна~$1$.}
 \end{itemize}

 Из условия следует, что $G/A$~--- $\pi'$-группа. Поэтому группа $$ N_G(S)A/A\cong N_G(S)/N_A(S)$$  и ее гомоморфный образ
  \begin{multline*}
 N_G(S)/N_A(S)C_G(S)\cong\\ (N_G(S)/C_G(S))/(N_A(S)C_G(S)/C_G(S))\cong\\ \Aut_G(S)/\Aut_A(S)
 \end{multline*}
  также являются $\pi'$-группами. Из определения множества $\Omega$ следует, что $\Aut_A(S)$ стабилизирует любой элемент из $\Omega$, поэтому длина любой орбиты на $\Omega$ группы $\Aut_G(S)$ делит $\pi'$-число
 $|\Aut_G(S):\Aut_A(S)|$ и следовательно сама является $\pi'$-чи\-слом. Если число является одновременно $\pi$- и $\pi'$-чи\-слом, то оно равно~$1$. Тем самым утверж\-де\-ние~$2^\circ$ доказано.

 \medskip

 Из $1^\circ$ и $2^\circ$ заключаем, что

 \begin{itemize}
   \item[$3^\circ$] {\it существует $H\in \Hall_{\mathfrak{X}}(S)$ с тем свойством, что для любого $\gamma\in\Aut_G(S)$ подгруппа $H^\gamma$ сопряжена в~$S$ с~$H$.}
 \end{itemize}

 Теперь мы можем утверждать следующее:
\begin{itemize}
   \item[$4^\circ$] {\it в любой минимальной нормальной подгруппе $M$ группы $G$ такой, что $M\le N$, существует $V_M\in\Hall_{\mathfrak{X}}(M)$, для которой $V_M^M=V_M^G$. }
 \end{itemize}

Не уменьшая общности, можно считать, что $M=\langle S^G\rangle$. Утверждение~$4^\circ$ следует из леммы~\ref{Constrain}.

\begin{itemize}
   \item[$5^\circ$] {\it   существует $V\in\Hall_{\mathfrak{X}}(N)$, для которой $V^N=V^G$. }
 \end{itemize}

Обозначим через $\Lambda$
множество всех содержащихся в
$N$  минимальных нормальных подгрупп группы $G$. Из условия следует, что
$$
N=\prod_{M\in \Lambda
} M,
$$
где произведение прямое. Для каждой такой $M$ выберем в соответствии с $4^\circ$ подгруппу $V_M\in\Hall_{\mathfrak{X}}(M)$ так, чтобы $V^M_M=V^G_M$. Пусть
$$
V=\langle V_M\mid M\in \Lambda \rangle.
$$
Тогда $V\in\Hall_{\mathfrak{X}}(N)$. Возьмем    $g\in G$. В любой $M\in\Lambda$ найдется элемент $x_M^{\vphantom{a}}$ такой, что $V_M^{g_{\vphantom{M}}^{\vphantom{a}}}=V_M^{x_M^{\vphantom{a}}}$. Положим $$x=\prod_{M\in \Lambda%\Xi%\Upsilon
} x_M^{\vphantom{a}}.
$$
Ясно, что $x\in N$ и $$V_M^{x_{\vphantom{M}}^{\vphantom{a}}}=V_M^{x_M^{\vphantom{a}}}=V_M^{g_{\vphantom{M}}^{\vphantom{a}}}$$ для любой $M\in\Lambda$. Значит,
$$
V^g=\langle V_M^g\mid M\in \Lambda \rangle=\langle V_M^x\mid M\in \Lambda \rangle=V^x.
$$
Утверждение $5^\circ$ доказано.

\begin{itemize}
   \item[$6^\circ$] {\it  Группы $N_A(V)$ и $N_G(V)$ являются ${\mathfrak{X}}$-отделимыми.}
 \end{itemize}
 Рассмотрим нормальный ряд $$N_G(V)\trianglerighteqslant N_A(V)\trianglerighteqslant N_N(V)\trianglerighteqslant V\trianglerighteqslant 1 $$ и его секции.
 Секция $N_G(N)/N_A(V)$ изоморфна подгруппе $N_G(V)A/A$ в~${\mathfrak{X}}'$-груп\-пе $G/A$ и значит сама является ${\mathfrak{X}}'$-груп\-пой. Аналогично $$N_A(V)/N_N(V)\cong N_A(V)N/N\le A/N=KN/N\cong K/(K\cap N),$$
 откуда $N_A(V)/N_N(V)\in{\mathfrak{X}}$. Так как $V\in\Hall_{\mathfrak{X}}(N)$, имеем $N_N(V)/V\in {\mathfrak{X}}'$. Наконец, $V\in{\mathfrak{X}}$. Таким образом, $6^\circ$ доказано.

\medskip

 Теперь из $5^\circ$ и $6^\circ$ выведем доказываемое предложение. Из $5^\circ$ вытекает, что $V^A=V^N$. По лемме~\ref{Lifting}  существует $L\in \Hall_{\mathfrak{X}}(A)$ такая, что $V=L\cap N$. Покажем, что $L$ удовлетворяет заключению предложения. Достаточно доказать включение $G\le AN_G(L)$. Ясно, что $L\le N_A(V)$, т.\,е.
  $L$~--- ${\mathfrak{X}}$-холлова подгруппа ${\mathfrak{X}}$-отделимой нормальной подгруппы $N_A(V)$ группы $N_G(V)$. Из сопряженности ${\mathfrak{X}}$-холловых подгрупп в ${\mathfrak{X}}$-отделимых группах получаем $$L^{N_G(V)}=L^{N_A(V)},\text{ откуда }N_G(V)\le N_A(V)N_G(L).$$ Теперь в силу $5^\circ$ имеем
$$
G=NN_G(V)\le NN_A(V)N_G(L)\le AN_G(L),
$$
что завершает доказательство предложения.
\qed\par\medskip

\noindent{\bf Замечание.} В идейном плане доказательство предложения~\ref{XFrattini} не содержит в себе ничего нового и с минимальными изменениями воспроизводит рассуждения из~\cite{Frattini}. По-видимому, можно доказать точный аналог основного результата~\cite{Frattini} для $\mathfrak{X}$-холловых подгрупп, а именно:
  \begin{itemize}
   \item[$\bullet$] {\it если $G\in\E_{\mathfrak{X}}$ и $A\trianglelefteqslant
   G$, то $G=AN_G(H)$ для некоторой $H\in\Hall_{\mathfrak{X}}(A)$.}
 \end{itemize}
  Точное следование рассуждениям цитированной работы потребовало использовать критерий существования $\mathfrak{X}$-холловых подгрупп, аналогичный полученному в~\cite{RV_Exist,Gross_Exist} и доказанный в~\cite[теорема~3.1]{Ved} (см. также   \cite[теорема~2]{BKh}).
 %\end{remark}
\section{О простых группах с девятью классами сопряженности ${\mathfrak{X}}$-холловых подгрупп}

\begin{propos}\label{NotMX}  {\it Пусть ${\mathfrak{X}}$~--- полный класс конечных групп, $S$~--- неабелева простая группа и $\h_{\mathfrak{X}}(S)=9$. Тогда $S\notin\M_{\mathfrak{X}}$.}\end{propos}

\noindent{\it Доказательство.}
Допустим, $S\in\M_{\mathfrak{X}}$. Пусть $\pi=\pi({\mathfrak{X}})$.
В соответствии с леммой~\ref{ClassNumberX=9} считаем, что
$$S= \PSp_{2n}(q)\cong \PSp(V)\quad\text{и}\quad\pi({\mathfrak{X}})\cap\pi(S)\subseteq\pi(\SL_2(q))=\pi(q^2-1).$$
Далее, любая ${\mathfrak{X}}$-холлова подгруппа группы $\PSp_{2n}(q)$ содержится в стабилизаторе~$M$ разложения ассоциированного с $\PSp_{2n}(q)$ пространства $V$ в ортогональную сумму
       $$V=V_1\perp \dots\perp V_n$$ невырожденных изометричных подпространств размерности~$2$. Существует подгруппа $A\trianglelefteqslant M$, такая, что $A=L_1\dots L_n$, где $L_i\cong \Sp(V_i)\cong \Sp_2(q)\cong \SL_2(q)$, $[L_i,L_j]=1$, $i,j=1,\dots,n$, $i\ne j$, и $M/A\cong \Sym_n$.

 Имеет
место один из двух случаев:
\begin{itemize}
  \item[$(1)$] $\pi({\mathfrak{X}})\cap\pi(S)=\pi({\mathfrak{X}})\cap\pi(\SL_2(q))=\{2,3\}$ и $n\in\{5,7\}$. При этом ${\mathfrak{X}}$-холловы подгруппы любой группы $L_i\cong\SL_2(q)$~--- это, в точности, обобщенные группы кватернионов порядка $48$ и группы вида $2.\Sym_4$.
  \item[$(2)$] $\pi({\mathfrak{X}})\cap\pi(S)=\pi({\mathfrak{X}})\cap\pi(\SL_2(q))=\{2,3,5\}$ и   $n=7$. При этом ${\mathfrak{X}}$-холловы подгруппы в любой $L_i\cong\SL_2(q)$~--- это, в точности, обобщенные группы кватернионов порядка $120$ и группы вида $2.\Alt_5$, а ${\mathfrak{X}}$-холлова подгруппа в $M/A\cong \Sym_7$ изоморфна~$\Sym_6$.
\end{itemize}
В каждом из случаев выберем в $S$ подгруппу $U$, как описано ниже.

Рассмотрим случай $(1)$. В группе $S$ имеется подгруппа вида $$\Sp_6(q)\circ \Sp_{2(n-3)}(q),$$ стабилизирующая в $S$ невырожденное подпространство размерности $6$ и его ортогональное дополнение, а значит есть подгруппа, изоморфная $\Sp_6(q)$. В~ней для любого $\varepsilon\in\{+,-\}$ есть подгруппа\footnote{Здесь и далее мы следуем стандартной для конечных классических групп практике (см.~\cite{Bray}, например), полагая $\GL_m^+(q)=\GL_m(q)$, $\SL_m^+(q)=\SL_m(q)$, $\GL_m^-(q)=\mathrm{GU}_m(q)$ и~$\SL_n^-(q)=\mathrm{SU}_n(q)$.}  $\GL^\varepsilon_3(q).2$, см.~\cite[таблица~8.28]{Bray}, причем $\varepsilon$ можно выбрать так, чтобы число $q-\varepsilon1$ делилось на~$3$. При таком выборе $\varepsilon$ в соответствии с~\cite[таблицы~8.3 и~8.5]{Bray} в подгруппе $\SL^\varepsilon_3(q)\le \GL^\varepsilon_3(q).2$ возьмем $\{2,3\}$-подгруппу $$U:= 3_+^{1+2}:Q_8.$$ В силу разрешимости $U\in{\mathfrak{X}}$. Так как $S\in\M_{\mathfrak{X}}$, имеем $U\le H$ для некоторой $H\in\Hall_{\mathfrak{X}}(S)$. Выберем подгруппу $M$ и в ней нормальную подгруппу $A$, как описано выше, с тем, чтобы $H\le M$. Пусть
$$\overline{\phantom{x}}:M\rightarrow M/A$$ обозначает канонический эпиморфизм. Тогда $$\overline{U}\le \overline{H}\le \overline{M}\cong \Sym_n,\text{ где }n\in\{5,7\}.$$
С другой стороны, $\overline{H}\cong H/(H\cap A)$ и $$\overline{U}\cong U(H\cap A)/(H\cap A)%\cong U/(U\cap H\cap A)
.$$  Выберем в $H\cap A$ характеристические подгруппы $B$, $C$ и $D$, определенные следующими равенствами:
$$
B:=O_2(H\cap A),\quad C/B:=O_3((H\cap A)/B)\text{ и } D/C:=O_2((H\cap A)/C).
$$
По определению $B\le C\le D$.
Как следует из свойств ${\mathfrak{X}}$-холловых подгрупп группы~$S$, подгруппа $H\cap A$ порождается попарно перестановочными ${\mathfrak{X}}$-хол\-ло\-вы\-ми подгруппами сомножителей $L_i$, каждая из которых либо является обобщенной  $\{2,3\}$-группой кватернионов, либо изоморфна $2.\Sym_4$. Поэтому ясно, что $D=H\cap A$ и следовательно, $\overline{U}\cong U/(U\cap D)$. Так как $O_2(U)=1$ имеем $$U\cap B=1\text{ и }U\cong UB/B.$$ Поскольку силовские $3$-подгруппы в каждом из сомножителей, образующих группу $H\cap A$, являются циклическими группами порядка $3$, силовская $3$-под\-груп\-па группы $H\cap A$, изоморфная $C/B$, абелева, и ее секция $$(UB/B)\cap (C/B)=(U\cap C)B/B\cong U\cap C$$ является нормальной абелевой 3-подгруппой группы $UB/B\cong U$. Отсюда следует, что она содержится в $Z(O_3(UB/B))$, так как $$UB/B\cong 3_+^{1+2}:Q_8,$$ а $Q_8$ действует неприводимо на факторгруппе группы $3_+^{1+2}$ по центру. Значит,  $$\text{либо }UC/C\cong U/(U\cap C)\cong 3_+^{1+2}:Q_8,\text{ либо }UC/C\cong 3^{2}:Q_8.$$ Наконец, из того, что $O_2(3_+^{1+2}:Q_8)=1$ и $O_2(3^{2}:Q_8)=1$ следует, что $$(UC/C)\cap D/C=1\text{ и }\overline{U}=UD/D\cong UC/C.$$ Но $\overline{U}$ (и следовательно ее подгруппа~$Q_8$) изоморфна подгруппе группы $\Sym_n$, для $n\in\{5,7\}$. Но довольно очевидно, что у группы $Q_8$ нет точных подстановочных представлений степени меньше~$8$. Противоречие.

Рассмотрим случай $(2)$. Рассуждая, как в случае $(1)$, в группе $S=\PSp_{14}(q)$ находим подгруппу, изоморфную~$\Sp_{10}(q)$. Поскольку $q^2-1$ делится на~$5$, возьмем $\varepsilon\in\{+,-\}$ так, чтобы $q-\varepsilon1$ делилось на~$5$. В группе $\Sp_{10}(q)$ a значит и в группе $S$  есть подгруппа~$\GL_5^\varepsilon(q).2$, см.~\cite[таблица~8.64]{Bray}, а~в~ней подгруппа~$\SL_5^\varepsilon(q)$. В~$\SL_5^\varepsilon(q)$, в свою очередь, найдется подгруппа $$U:=5_+^{1+2}.\Sp_2(5),$$ см.~\cite[таблицы~8.18 и~8.20]{Bray}, причем, так как $\Sp_2(5)\cong 2.\Alt_5$, имеем $U\in{\mathfrak{X}}$. Так как $S\in\M_{\mathfrak{X}}$, имеем $U\le H$ для некоторой $H\in\Hall_{\mathfrak{X}}(S)$. Выберем подгруппу $M$ и в ней нормальную подгруппу $A$, как описано выше, с тем, чтобы $H\le M$. Пусть
$$\overline{\phantom{x}}:M\rightarrow M/A$$ обозначает канонический эпиморфизм. Тогда $$\overline{U}\le \overline{H}\le \overline{M}\cong \Sym_7.$$ Следовательно, $|\overline{U}|_5\le 5$. Из строения группы $U$ видно, что любой гомоморфный образ группы $U$, порядок которого не делится на $5^2$, является образом группы $\Sp_2(5)\cong \SL_2(5)$. Поэтому экстраспециальная подгруппа $5_+^{1+2}$ группы $U$ должна лежать ядре гомоморфизма $\overline{\phantom{x}}$, а значит, лежать в $U\cap A$. Но силовская $5$-подгруппа группы $A$ абелева (в каждом сомножителе $L_i$ силовская $5$-подгруппа имеет порядок~5). Противоречие.
\qed\par\medskip

\section{Доказательство теоремы~\ref{Reduktionssatz_main}  и ее следствий}

\noindent{\it Доказательство теоремы~}\ref{Reduktionssatz_main}. Утверждение ``{\it если}'' теоремы~\ref{Reduktionssatz_main} доказано в~\cite[теорема~1]{GRV_Reducktionssatz} (см.~лемму~\ref{RedSatz}). Предположим, утверждение ``{\it только если}'' неверно. Тогда существует группа~$G$ со следующими свойствами:

\begin{itemize}
  \item[$(a)$]\label{a} в $G$ имеется нормальная подгруппа $N$ такая, что ${\mathrm{k}}_{\mathfrak{X}}(N)>1$, но для пары $(G,N)$ выполнена редукционная ${\mathfrak{X}}$-теорема, т.\,е. ${\mathrm{k}}_{\mathfrak{X}}(G)={\mathrm{k}}_{\mathfrak{X}}(\overline{G})$, где черта $$\overline{\phantom{x}}:G\rightarrow G/N$$ обозначает канонический эпиморфизм;
  \item[$(b)$]\label{b} $G$ имеет наименьший порядок среди групп, обладающих свойством~$(a)$.
  %\item[$(c)$]
\end{itemize}

Напомним, что редукционная ${\mathfrak{X}}$-теорема для пары $(G,N)$ влечет следующие свойства:
 {\it
 \begin{itemize}
   \item[$1^\circ$]\label{1circ} если  $K\in\m_{\mathfrak{X}}(G)$, то $\overline{K}\in\m_{\mathfrak{X}}(\overline{G})$, и
   \item[$2^\circ$]\label{2circ} если для $K,L\in\m_{\mathfrak{X}}(G)$ подгруппы $\overline{K}$ и $\overline{L}$ сопряжены в $\overline{G}$ $($например, совпадают$)$, то $K$ и $L$ сопряжены в~$G$.
 \end{itemize}}

Заметим, что если ${M\ne 1}$~--- нормальная подгруппа в~$G$, причем ${M\leqslant N}$, то, поскольку $G/N$ является гомоморфным образом группы $G/M$, имеем
 $$
 {\mathrm{k}}_{\mathfrak{X}}(G)={\mathrm{k}}_{\mathfrak{X}}(G/N)\le{\mathrm{k}}_{\mathfrak{X}}(G/M)\le{\mathrm{k}}_{\mathfrak{X}}(G).
 $$
 Отсюда следует редукционная ${\mathfrak{X}}$-теорема для пар $(G/M,N/M)$ и $(G,M)$. В~силу свойства $(b)$ и того, что ${|G/M|<|G|}$, имеем ${{\mathrm{k}}_{\mathfrak{X}}(N/M)=1}$. Поэтому если ${{\mathrm{k}}_{\mathfrak{X}}(M)=1}$, то из леммы~\ref{RedSatz} следует, что $${\mathrm{k}}_{\mathfrak{X}}(N)={\mathrm{k}}_{\mathfrak{X}}(N/M)=1.$$
 Следовательно, мы можем считать, что
  {\it  \begin{itemize}
   \item[$3^\circ$]\label{3circ} $N$~---  минимальная нормальная подгруппа в~$G$, причем $N$ неабелева,
 \end{itemize}}
\noindent так как ${\mathrm{k}}_{\mathfrak{X}}(N)>1$. Значит,
 $$
 N=S_1\times\dots\times S_n
 $$
 для некоторых неабелевых простых подгрупп $S_1,\dots ,S_n$, сопряженных в~$G$. Пусть $S$~--- одна из них.

 Мы получим противоречие, изучая действие группы $\Aut_G(S)$ на множестве $$\Delta:=\Hall_{\mathfrak{X}}(S)/S$$ классов сопряженности ${\mathfrak{X}}$-холловых подгрупп группы~$S$.
 Как следует из леммы~\ref{ClassNumberX},
 {\it
 \begin{itemize}
   \item[$4^\circ$] $
|\Delta|= \h_{\mathfrak{X}}(S)\in\{0,1,2,3,4,9\}.
 $
 \end{itemize}}
\noindent  Исключим все шесть возможностей. Сначала установим, что
  {\it  \begin{itemize}
   \item[$5^\circ$] $\h_{\mathfrak{X}}(S)\ne 0$.
 \end{itemize}}
\noindent  Для этого докажем, что
 \begin{itemize}
   \item[$6^\circ$] {\it если $K\in\m_{\mathfrak{X}}(G)$, то $K\cap N\in\Hall_{\mathfrak{X}}(N)$ и $K\in\Hall_{\mathfrak{X}}(KN)$.}
 \end{itemize}
\noindent Отсюда будет следовать, что $\Hall_{\mathfrak{X}}(N)\ne\varnothing$ и $$\varnothing\ne\{H\cap S\mid H\in\Hall_{\mathfrak{X}}(N)\}\subseteq\Hall_{\mathfrak{X}}(S),$$ т.\,е. ${\h_{\mathfrak{X}}(S)\ne 0}$. Заметим, что  из $K\cap N\in\Hall_{\mathfrak{X}}(N)$ следует, что  $K\in\Hall_{\mathfrak{X}}(KN),$ поскольку
$$
|KN:K|=\frac{|K||N|}{|K\cap N|}:|K|=|N:(K\cap N)|.
$$

 Выберем произвольно $K\in\m_{\mathfrak{X}}(G)$, $p\in\pi({\mathfrak{X}})$ и $P\in\Syl_p(N)$. Тогда $P\in{\mathfrak{X}}$. Достаточно установить, что $P$ сопряжена с подгруппой из $K$. Положим $$A:=KN.$$ Из аргумента Фраттини~\cite[гл.~A,~(6,3)]{DH} следует, что ${A={ N}_{A}(P)N}$. Поэтому $$\overline{K}=\overline{A}=\overline{{ N}_{A}(P)}$$ и согласно лемме~\ref{FrattiniSubgr} имеем $\overline{{ N}_{A}(P)}=\overline{U}$ для некоторой $U\in\m_{\mathfrak{X}}({{ N}_{A}(P)})$. Так как $U$ нормализует ${\mathfrak{X}}$-подгруппу $P$, имеем $P\le U$. Теперь погрузим $U$ в максимальную ${\mathfrak{X}}$-подгруппу $L$ группы~$G$. В силу $1^\circ$ имеем $\overline{L}\in\m_{\mathfrak{X}}(\overline{G})$. Кроме того,
$$
\overline{K}=\overline{{ N}_{A}(P)}=\overline{U}\le \overline{L}.
$$
 Аналогично в силу $1^\circ$ и того, что $K\in\m_{\mathfrak{X}}(G)$, имеем $\overline{K}\in\m_{\mathfrak{X}}(\overline{G})$. Теперь  $\overline{K}=\overline{L}$ и~ по~$2^\circ$ заключаем, что $L^g=K$ для некоторого $g\in G$. Поэтому $$P^g\le U^g\le L^g=K,$$ как и утверждалось.

 \medskip

 Отметим ряд других следствий утверждения~$6^\circ$. Мы утверждаем, что
 \begin{itemize}
   \item[$7^\circ$] {\it всякая ${\mathfrak{X}}$-подгруппа группы $G$ нормализует элемент из $\Hall_{\mathfrak{X}}(N)$; в частности}
   \item[$8^\circ$] {\it $N\in\M_{\mathfrak{X}}$ и $S\in\M_{\mathfrak{X}}$;}
   \item[$9^\circ$] {\it всякая ${\mathfrak{X}}$-подгруппа группы $\Aut_G(S)$ стабилизирует некоторый элемент из~$\Delta$.}
   \end{itemize}

\noindent   Докажем  $7^\circ$ и $8^\circ$. Если $U$~--- ${\mathfrak{X}}$-подгруппа группы $G$, то $U\le K$ для некоторой $K\in\m_{\mathfrak{X}}(G)$ и $U$ нормализует $K\cap N\in\Hall_{\mathfrak{X}}(N)$ согласно~$6^\circ$. Утверждение $7^\circ$ доказано. Если при этом взять $U\in\m_{\mathfrak{X}}(N)$, то $U(K\cap N)$~--- ${\mathfrak{X}}$-подгруппа в~$N$, содержащая $U\in \m_{\mathfrak{X}}(N)$ и $K\cap N\in\Hall_{\mathfrak{X}}(N)\subseteq \m_{\mathfrak{X}}(N)$, поэтому
$$
U=U(K\cap N)=K\cap N\in\Hall_{\mathfrak{X}}(N).
$$
Тем самым доказано, что ${\m_{\mathfrak{X}}(N)=\Hall_{\mathfrak{X}}(N)}$. Значит ${N\in\M_{\mathfrak{X}}}$. Как следует из леммы~\ref{HallSubgroup} нормальная подгруппа $\M_{\mathfrak{X}}$-группы является $\M_{\mathfrak{X}}$-группой. Отсюда получаем $S\in\M_{\mathfrak{X}}$, как и утверждается в~$8^\circ$.

Докажем $9^\circ$. Обозначим через $$\alpha: N_G(S)\rightarrow \Aut_G(S)$$ эпиморфизм, который сопоставляет элементу $g\in N_G(S)$ автоморфизм группы~$S$, действующий по правилу $x\mapsto x^g$ для всех ${x\in S}$. По лемме~\ref{FrattiniSubgr} произвольная ${\mathfrak{X}}$-подгруппа $T$ группы $\Aut_G(S)$ имеет вид~$U^\alpha$, где $U$~--- некоторая ${\mathfrak{X}}$-подгруппа группы~$N_G(S)$. Возьмем произвольно $t\in T$. Тогда найдется  $g\in U$ с тем свойством, что $$x^t=x^{g^\alpha}=x^g$$ для всех $x\in S$. В силу~$7^\circ$ подгруппа $U$ нормализует некоторую $H\in\Hall_{\mathfrak{X}}(N)$. Так как $U$ нормализует также подгруппу $S$, имеем
$$
(H\cap S)^{t}=(H\cap S)^{g}= H\cap S.
$$ Следовательно, $T$ оставляет инвариантным класс сопряженности ${\mathfrak{X}}$-холловых под\-групп группы $S$, содержащий $H\cap S\in\Hall_{\mathfrak{X}}(S).$

  {\it  \begin{itemize}
   \item[$10^\circ$] $\h_{\mathfrak{X}}(S)\ne 1$.
 \end{itemize}}
\noindent

Действительно, если $\h_{\mathfrak{X}}(S)=1$, то $S\in\C_{\mathfrak{X}}$ и в силу $8^\circ$ получаем $$S\in\C_{\mathfrak{X}}\cap \M_{\mathfrak{X}}=\D_{\mathfrak{X}}.$$ Но тогда $${\mathrm{k}}_{\mathfrak{X}}(S_1)=\ldots={\mathrm{k}}_{\mathfrak{X}}(S_n)=1$$ и ${\mathrm{k}}_{\mathfrak{X}}(N)=1$ по лемме~\ref{RedSatz}, вопреки~$(a)$.

 {\it  \begin{itemize}
   \item[$11^\circ$] $\h_{\mathfrak{X}}(S)\ne 9$,
 \end{itemize}}
 \noindent  поскольку $S\in\M_{\mathfrak{X}}$, а $S\notin\M_{\mathfrak{X}}$ при $\h_{\mathfrak{X}}(S)=9$ по предложению~\ref{NotMX}.

\begin{itemize}
   \item[$12^\circ$] {\it $N_G({KN})\in\C_{\mathfrak{X}}$  и $KN\in\C_{\mathfrak{X}}$ для любой $K\in\m_{\mathfrak{X}}(G)$;}
   \end{itemize}
    Пусть $K\in\m_{\mathfrak{X}}(G)$. Положим $$A:=KN\text{ и }B:=N_G(A).$$ Из того, что $\overline{A}=\overline{K}\in \m_{\mathfrak{X}}(\overline{G})$, заключаем, что $B/A\cong N_{\overline{G}}(\overline{K})/\overline{K}$~--- $\pi'$-группа. Сле\-до\-ва\-тель\-но, $$\Hall_{\mathfrak{X}}(A)=\Hall_{\mathfrak{X}}(B).$$ В соответствии с~$6^\circ$ имеем ${K\in\Hall_{\mathfrak{X}}(A)}$. В частности, ${\Hall_{\mathfrak{X}}(A)\ne\varnothing}$.  Значит, $$A,B\in\E_{\mathfrak{X}}\text{ и }\h_{\mathfrak{X}}(A)\ge\h_{\mathfrak{X}}(B)\ge 1.$$ Допустим, $L\in\Hall_{\mathfrak{X}}(A)$ и покажем, что $L$ и $K$ сопряжены в~$B$. Для начала $A=KN=LN$ и значит ${\overline{K}=\overline{L}}$. Поэтому $$\overline{L}\in\m_{\mathfrak{X}}(\overline{G})\text{ и }L\cap N\in\Hall_{\mathfrak{X}}(N)\subseteq\m_{\mathfrak{X}}(N).$$ Отсюда $L\in\m_{\mathfrak{X}}(G)$. В соответствии с~$2^\circ$,  равенство $\overline{K}=\overline{L}$ влечет сопряженность подгрупп $K$ и $L$ в~$G$ и следовательно в $B=N_G(KN)=N_G(KL)$. Мы доказали, что группа $B$ транзитивно действует сопряжениями на непустом множестве $\Hall_{\mathfrak{X}}(A)=\Hall_{\mathfrak{X}}(B)$. Значит, $N_G(KN)=B\in \C_{\mathfrak{X}}$.

    Из предложения~\ref{XFrattini} следует, что {\it существует} $L\in\Hall_{\mathfrak{X}}(A)$, для которой выполнено равенство $B=N_B(L)A$. Теперь, поскольку $\Hall_{\mathfrak{X}}(A)=\Hall_{\mathfrak{X}}(B)$ и $B\in\C_{\mathfrak{X}}$, это равенство выполнено {\it для любой} $L\in\Hall_{\mathfrak{X}}(A)$. Из $B\in \C_{\mathfrak{X}}$ следует также, что для любой $L\in\Hall_{\mathfrak{X}}(A)$ существует элемент $b\in B$ такой, что $K=L^b$. При этом $b=ga$ для некоторых $g\in N_B(L)$ и $a\in A$. Значит,
    $$
    K=L^b=L^{ga}=L^a,\text{ где } a\in A,
    $$
    и тем самым $KN=A\in\C_{\mathfrak{X}}$, как и утверждалось.

\medskip

Утверждение~$12^\circ$ позволяет установить следующий факт, который оказывается решающим в исключении оставшихся случаев:

     \begin{itemize}
   \item[$13^\circ$] {\it  для любой ${K\in\m_{\mathfrak{X}}(G)}$ группа $\Aut_K(S)$ стабилизирует ровно один элемент из~$\Delta$.}
 \end{itemize}

 Возьмем $K\in\m_{\mathfrak{X}}(G)$. Существование в $\Delta$ неподвижной точки относительно $\Aut_K(S)$ следует из~$9^\circ$. Пусть $H\in\Hall_{\mathfrak{X}}(S)$. Мы докажем $13^\circ$, если установим, что из инвариантности класса $H^S\in\Hall_{\mathfrak{X}}(S)/S$ относительно $\Aut_K(S)$ следует, что $H\in(K\cap S)^S$. Пусть, как и раньше, $A=KN$. По лемме \ref{AutK=AutNK}  имеем $${H\vphantom{\left(H^S\right)}}^{\Aut_A(S)}=\left(H^S\right)^{\Aut_K(S)}=H^S.$$
  %Заметим, что $N\le N_X(S)$, поэтому $N_X(S)=NN_K(S)$, причем подгруппа $N$ стабилизирует каждый класс сопряженности подгрупп в~$S$. Поэтому
% $$H^{\Aut_X(S)}=H^{N_X(S)}=H^{NN_K(S)}=\left(H^S\right)^{N_K(S)}=\left(H^S\right)^{\Aut_K(S)}=H^S.$$
 Теперь из леммы~\ref{Constrain} следует, что для $$M:=\langle S^A\rangle=\langle S^K\rangle$$ существует $L\in\Hall_{\mathfrak{X}}(KM)$ такая, что $H=L\cap S$. При этом $|L|=|K|$ и значит $L\in\Hall_{\mathfrak{X}}(A)$. Но $A=KN\in\C_{\mathfrak{X}}$ в силу~$12^\circ$, поэтому существуют $u\in K$ и $v\in N$ такие, что
 $$
 L=K^{uv}=K^v.
 $$
 При этом понятно, что если $w\in S$~--- проекция элемента $v$ на $S$,  (напомним, $S$~--- один из прямых сомножителей $S_1,\dots,S_n$, произведение которых дает группу~$N$), то $$H=L\cap S=K^v\cap S=(K\cap S)^v=(K\cap S)^{w}\in (K\cap S)^{S}.$$ Тем самым $13^\circ$ доказано.

 \medskip

 Из $4^\circ$, $5^\circ$, $10^\circ$ и $11^\circ$ следует, что $\h_{\mathfrak{X}}(S)\in\{2,3,4\}$. Покажем, что

 {\it  \begin{itemize}
   \item[$14^\circ$] $\h_{\mathfrak{X}}(S)\ne 2$.
 \end{itemize}}
\noindent  В самом деле, зафиксируем некоторую $K\in\m_{\mathfrak{X}}(G)$. Согласно $13^\circ$ группа $\Aut_K(S)$ стабилизирует ровно один элемент из $\Delta$. Но при  $\h_{\mathfrak{X}}(S)=2$ группа $\Aut_K(S)$ ста\-би\-ли\-зи\-рует и оставшийся элемент. Противоречие.

С учетом сказанного и в силу леммы~\ref{ClassNumberX} мы можем считать, что
{\it  \begin{itemize}
   \item[$15^\circ$] $\h_{\mathfrak{X}}(S)\in\{3,4\}$ и $2,3\in\pi({\mathfrak{X}})$.
 \end{itemize}}
 Теперь по теореме Бернсайда \cite[гл.~I, разд.~2]{DH} и в силу того, что всякая разрешимая $\pi({\mathfrak{X}})$-группа является ${\mathfrak{X}}$-группой,
 {\it  \begin{itemize}
   \item[$16^\circ$] всякая $\{2,3\}$-группа является ${\mathfrak{X}}$-группой.
 \end{itemize}}

 Мы получим окончательное противоречие с $9^\circ$,  доказывающее теорему, установив, что
 {\it  \begin{itemize}
   \item[$17^\circ$] группа $\Aut_G(S)$ содержит ${\mathfrak{X}}$-подгруппу, транзитивно действующую на~$\Delta$.
 \end{itemize}}

\noindent   Действие группы $\Aut_G(S)$ на множестве $\Delta$ индуцирует гомоморфизм
 $$
 {}^*:\Aut_G(S)\rightarrow \Sym(\Delta),\text{ где }\Sym(\Delta)\cong
 \left\{\begin{array}{cr}
                                              \Sym_3, & \text{ если }\quad \h_{\mathfrak{X}}(S)=3, \\
                                               \Sym_4, & \text{ если }\quad \h_{\mathfrak{X}}(S)=4.
                                            \end{array}
                                            \right.
 $$
 Возьмем произвольный класс $$H^S=\{H^x\mid x\in S\}\in\Delta$$ с представителем ${H\in\Hall_{\mathfrak{X}}(S)}$, и пусть $H\le K$ для некоторой $K\in\m_{\mathfrak{X}}(G)$. Поскольку подгруппа $H=K\cap S$ инвариантна относительно $N_K(S)$, класс сопряженности $H^S$ стабилизируется группой $\Aut_K(S)$ и согласно~$11^\circ$ эта группа действует без не\-под\-виж\-ных точек на множестве $$\Gamma:=\Delta\setminus\{H^S\}$$ оставшихся двух или трех классов. Из того, что $$|\Gamma|=\h_{\mathfrak{X}}(S)-1=\left\{\begin{array}{cr}
                                              2, & \text{ если }\quad \h_{\mathfrak{X}}(S)=3, \\
                                               3, & \text{ если }\quad \h_{\mathfrak{X}}(S)=4,
                                            \end{array}
                                            \right.$$ видно, что действие $\Aut_K(S)$ на $\Gamma$ должно быть транзитивным. Но тогда и действие стабилизатора в $\Aut_G(S)$
точки $H^S$ на $\Gamma$ также транзитивно, поскольку этот стабилизатор содержит подгруппу~$\Aut_K(S)$. Следовательно, $\Aut_G(S)^*$~--- транзитивная и даже дважды транзитивная подгруппа в $\Sym(\Delta)$, т.\,е.
$$
\Aut_G(S)^*\cong \left\{\begin{array}{rl}
                                              \Sym_3, & \text{ если }\quad \h_{\mathfrak{X}}(S)=3, \\
                                               \Alt_4\text{ или }\Sym_4, & \text{ если }\quad \h_{\mathfrak{X}}(S)=4.
                                            \end{array}
                                            \right.
$$
В любом случае $\Aut_G(S)^*$~--- $\{2,3\}$-группа и, следовательно, ${\mathfrak{X}}$-группа. По лемме~\ref{FrattiniSubgr} существует ${\mathfrak{X}}$-подгруппа $U$ группы $\Aut_G(S)$ такая, что $U^*=\Aut_G(S)^*$. Такая подгруппа $U$ транзитивно действует на $\Delta$, вопреки~$9^\circ$.

Теорема полностью доказана.\qed
\par\medskip

\noindent{\it Доказательство следствий}~\ref{Description}--\ref{Radical}. См. введение~\ref{Motivation},~\ref{Corollaries}.
\qed
\par\medskip

\noindent{\it Доказательство следствия}~\ref{Overgroups}. Пусть $H\in \m_{\mathfrak{X}}(G)$ и $H\leqslant K\leqslant G$. Пусть подгруппа $N\trianglelefteqslant G$ такова, что  $\mathrm{k}_{\mathfrak{X}}(G/N)=\mathrm{k}_{\mathfrak{X}}(G).$ Тогда $N\in\D_{\mathfrak{X}}$ по теореме~\ref{Reduktionssatz_main}. Покажем, что $$\mathrm{k}_{\mathfrak{X}}(K/(K\cap N))=\mathrm{k}_{\mathfrak{X}}(G).$$

Допустим, $N$ минимальна в $G$ как нормальная подгруппа. Тогда поскольку $N$ является прямым произведением своих нормальных подгрупп, сопряженных в $G$ и изоморфных некоторой простой группе $S$, согласно  \cite[теорема~2, лемма~2.28]{GRV_Reducktionssatz} либо $N\in\mathfrak{X}$, либо $N\in \D_\pi$ для $\pi=\pi(\mathfrak{X})$ и любая $\pi$-холлова подгруппа группы $N$ разрешима, в частности, принадлежит~$\mathfrak{X}$. В первом случае доказываемое очевидно. Во втором случае из леммы~\ref{HallSubgroup} следует, что $$H\cap N\in\Hall_{\mathfrak{X}}(N)\subseteq\Hall_\pi(N)\quad\text{и}\quad H\cap N\leqslant K\cap N\leqslant N.$$ Теперь $K\cap N\in\D_\pi$ согласно \cite[теорема~1.4]{VMR} и, так как $\pi$-холлова подгруппа группы $N$ принадлежит~$\mathfrak{X}$, имеем~$K\cap N\in \D_{\mathfrak{X}}$. Теперь $\mathrm{k}_{\mathfrak{X}}(K/(K\cap N))=\mathrm{k}_{\mathfrak{X}}(G)$ по теореме~\ref{Reduktionssatz_main}.

Общий случай теперь выводится из рассмотренного индукцией по~$|N|$ после перехода к  факторгруппе группы $G$ по минимальной нормальной подгруппе, содержащейся в~$N$, и последующего применения утверждения ``если'' в теореме~\ref{Reduktionssatz_main}.
\qed
\par\medskip

\noindent{\it Доказательство следствия}~\ref{IsoschematismsCategory}. Из определения отношения $\underset{{\mathfrak{X}}}{\equiv}$ следует, что $G_1\underset{{\mathfrak{X}}}{\equiv} G_2$ тогда и только тогда, когда полные $\mathfrak{X}$-редукции групп $G_1$ и $G_2$ изоморфны. Теперь все утверждения следствия~\ref{IsoschematismsCategory} очевидны.
\qed
\par\medskip

%\end{fulltext}

\medskip

\noindent {\bf Wenbin~Guo}\\
{\small 1. School of Science, Hainan University,
Haikou,} \\
{\small Hainan, 570228, P.R. China, and}\\
{\small 2. Department of Mathematics,}\\
{\small  University of Science and
Technology of China,}\\
{\small Hefei 230026, P. R. China}\\
{\small email: wbguo@ustc.edu.cn}\\ %или \email[r_rrrr@rrr.ru]{r\_rrrr@rrr.ru}, если в адресе есть нижнее подчеркивание
 ~\\
{\bf Danila O.~Revin}\\
{\small 1. Sobolev Institute of Mathematics,}\\
{\small  4, Koptyug av., Novosibirsk 630090}\\
{\small 2. Novosibirsk State University,}\\
{\small  1, Pirogova st. , Novosibirsk 630090}\\
{\small 3. Krasovsii Institute of Mathematics and Mechanics,}\\
{\small  16, S.~Kovalevskoy st., Yekaterinburg 620990, Russia}\\
{\small  email: revin@math.nsc.ru}
\end{document}